\documentclass[english,11pt]{article}
\usepackage{algpseudocode} 
\usepackage[ruled,linesnumbered]{algorithm2e}
\usepackage[english]{babel}
\usepackage{textcomp}
\usepackage{eurosym}
\usepackage{amsthm}
\usepackage[intlimits]{amsmath}
\usepackage{amssymb}
\usepackage{hyperref}
\usepackage{multirow}
\usepackage[table]{xcolor}
\usepackage{tabularx}

\usepackage{amsfonts}
\usepackage{paralist}
\usepackage{graphicx}
\usepackage{fancyhdr}
\usepackage{ulem}
\usepackage{epsfig}
\usepackage{color}

\newtheorem{thm}{Theorem}%[section]

\newtheorem{conj}[thm]{Conjecture}
\newtheorem{lem}[thm]{Lemma}
\newtheorem{tab}[thm]{Table}

\newtheorem{prop}[thm]{Proposition}
\newtheorem{exmp}[thm]{Example}

\newtheorem{defn}[thm]{Definition}
\theoremstyle{plain}

\newtheorem{rem}[thm]{Remark}

\theoremstyle{definition}

\numberwithin{equation}{section}

\theoremstyle{definition}

\newcommand{\EE}{\mathbb{E}}

\newcommand{\NN}{\mathbb{N}}

\newcommand{\PP}{\mathbb{P}}

\newcommand{\RR}{\mathbb{R}}

\newcommand{\TT}{\mathbb{T}}

\newcommand{\bB}{\mathcal{B}}
\newcommand{\cC}{\mathcal{C}}

\newcommand{\fF}{\mathcal{F}}

\newcommand{\hH}{\mathcal{H}}
\newcommand{\iI}{\mathcal{I}}

\newcommand{\mM}{\mathcal{M}}
\newcommand{\nN}{\mathcal{N}}

\newcommand{\pP}{\mathcal{P}}

\newcommand{\al}{\alpha}

\newcommand{\e}{\varepsilon}

\newcommand{\si}{\sigma}
\newcommand{\om}{\omega}
\newcommand{\Om}{\Omega}

\newcommand{\vp}{\varphi}

\newcommand{\ra}{\rightarrow}

\newcommand{\ti}{\tilde}

\newcommand{\id}{\ensuremath{id}}
\newcommand{\supp}{\ensuremath{supp}}
\newcommand{\lqq}{\leqslant}
\newcommand{\gqq}{\geqslant}

\DeclareMathOperator{\per}{Sym}

\usepackage[]{color}
\definecolor{DarkGreen}{rgb}{0.1,0.7,0.3}   %define a custom color
\definecolor{Yellow}{rgb}{0.7,0.7,0.7}   %define a custom color

\begin{document}
\thispagestyle{empty}

\title{Stochastic n-point D-bifurcations 
of stochastic L\'evy flows 
and their complexity on finite spaces}

\author{
Paulo Henrique da Costa\footnote{Universidade de Bras\'ilia, Bras\'ilia, Brazil; phcosta@unb.br} 
\hspace{1cm}
Michael A. H\"ogele\footnote{Universidad de los Andes, Bogot\'a, Colombia. ma.hoegele@uniandes.edu.co} 
\hspace{1cm}
Paulo R. Ruffino\footnote{State University of Campinas - UNICAMP, Brazil. ruffino@unicamp.br}
}

\maketitle
\begin{center}
At the occasion of Bj\"orn Schmalfu\ss's 65th birthday.\\[1cm]
\end{center}

\begin{abstract}

This article refines the classical notion of a stochastic D-bifurcation 
to the respective family of n-point motions for homogeneous Markovian stochastic semiflows, 
such as stochastic Brownian flows of homeomorphisms, and their generalizations. 
This notion essentially detects at which level %$k\lqq n$ 
the support of the invariant measure of the k-point bifurcation has more than one 
connected component. Stochastic Brownian flows and 
their invariant measures were shown by Kunita (1990) to be rigid, 
in the sense of being uniquely determined by the $1$-and $2$-point motions. Hence 
only stochastic n-point bifurcation of level $n=1$ or $n=2$ can occur. 
For general homogeneous stochastic Markov semiflows this turns out to be false.
This article constructs minimal examples of where this rigidity is false in general on finite space and studies the complexity of the resulting n-point bifurcations. 
\end{abstract}



% \tableofcontents
    
\section{Motivation}

This article refines the classical notion of a \textit{stochastic dynamical bifurcation} or \textit{stochastic D-bifurcation} 
given by L.~Arnold in \cite{Arnold} to the respective n-point motion for homogeneous Markovian stochastic semiflows of (continuous) mappings, such as stochastic Brownian flows of homeomorphisms,  
and their generalizations. The generalizations include coalescing stochastic flows, which go back to the PhD thesis by Arriata \cite{Arriata79}, but are an active field of research \cite{BGS15, DRS20, Do05, DK17, DV18, Ha84, LeJan-Raimond,Ri18}.
For this purpose we introduce the notion of a \textit{stochastic n-point D-bifurcation} 
in this setting. It essentially detects the lowest level $n$ 
at which the respective n-point motion of such a stochastic flow 
undergoes a stochastic D-bifurcation in the sense of 
a qualitative change of the support of the invariant measure. 
This concept is finer than the classical stochastic D-bifurcation, since trivially any stochastic D-bifurcation is equivalent to a stochastic $1$-point D-bifurcation. 
To the best of our knowledge, the concept of a stochastic n-point D-bifurcation 
is new in the literature and in our view it has the potential to provide 
interesting new insights into the tipping behavior of stochastic L\'evy (semi-)flows 
and more general systems given by L\'evy driven stochastic (partial) differential equations. 

It is known for a long time that the stochastic dynamics 
of a stochastic differential equation is only understood partially by the 
respective Markov semigroups and their generators.  
Indeed, these objects do not take into account the effect of the 
dependence structure between the trajectories of $n$ ensemble members 
with different initial conditions, that is, the respective n-point motion. 
This conceptual problem was solved with the introduction of 
the notion of the associated stochastic flow in the works  
by Elworthy, Baxendale, Bismut, Ikeda, Kunita, Watanabe among others 
(see e.g. \cite{Baxendale84, Bismut80, El78, IW89, Ku90, Ku19} and the references therein).  
Therefore the motivation for the generalization of stochastic D-bifurcation to n-point motions is twofold: 
\begin{enumerate}
 \item[I)] The 1- and 2-point rigidity of the laws of stochastic Brownian flows of homeomorphisms. 
 \item[II)] Discretization procedures for homogeneous Markovian semiflows of continuous functions.  
\end{enumerate}
Motivation I) is rooted in a sequence of results 
for homogeneous stochastic Brownian flows of homeomorphisms in Kunita \cite{Ku90}. 
The distribution of such a stochastic Brownian flow 
is uniquely determined by the laws of the families of the corresponding $1$-point and 
$2$-point motions, see Theorem 4.2.5 and formula (19) in Kunita \cite{Ku90}.  
This sort of 1- and 2-point rigidity of the law of the flow 
is due to the Gaussian nature of the marginal laws and their n-point extensions 
and can be read off directly from the structure of the infinitesimal generators. 
This rigidity carries over to any invariant distribution $\Pi$ of the Brownian flow 
(in the sense that the flow is $\Pi$-preserving) as follows: the flow is $\Pi$-preserving if and only if $\Pi$ is the invariant measure of the respective 1-point motion and $\Pi\otimes \Pi$ is the invariant measure of the 2-point motion (see Theorem 4.3.2(v) in \cite{Ku90}). In other words, the Gaussianity of Brownian flows imposes 
that the complete dependence structure of the n-point motion 
of the flow is uniquely determined by the respective infinitesimal 
covariances of 2-point motion contained in the n-point motion, 
that is to say, its 1 and 2-point characteristics. 
However, the 2-point characteristics can change the law of the flow, as shown by Baxendale in \cite{Baxendale91}, who studies the ergodicity of the 1- and 2-point Brownian motion on a torus.  
Homogeneous Markov semiflows of homeomorphisms, which we denote for convenience as stochastic L\'evy semiflows, generalize the notion of the respective stochastic Brownian flows by dropping the continuity assumption in time, cf. \cite{AK00, Fu91, FK85, Ku04}. 
However, in general, neither their (non-Gaussian) laws 
nor the respective invariant measures 
can be expected to be rigid in the sense of Brownian flows. 
This is due to the lack of continuity resulting in 
the non-local nature of the infinitesimal generator. 
Therefore, the law and the invariant measures of the respective 
n-point motion for $n>2$ of the flow provide new and finer information about the law 
and the invariant measure of the flow. 
It is therefore natural to ask for examples of stochastic L\'evy flows 
of homeomorphisms whose laws are underdetermined by their 2-point motions 
and which are minimal in some sense. 
While non-trivial stochastic 
Brownian flows are confined to spaces 
where Gaussian laws can be defined properly such as Hilbert spaces, 
% that is, essentially Hilbert spaces 
stochastic L\'evy flows exhibit jumps and 
have a rich behavior already on spaces of finite sets. 
In general, L\'evy driven stochastic differential 
equations yield stochastic L\'evy flows under rather 
restrictive conditions on the coefficients \cite{Applebaum, AK00, Ku04}. 
However, their easiest representatives are given 
as continuous-time homogeneous Markov chains with values 
in a finite state space. In order to 
settle ideas we focus in this article 
on finite state spaces. 

Motivation II) is a bit more far-flung. There is an ever growing necessity to detect (stochastic) bifurcations or tipping of high dimensional stochastic semiflows or even more general systems, such as stochastically perturbed general circulation models in climatology (see for instance \cite{IM02, Im01}). The notion of a stochastic n-point D-bifurcation over finite points seems to be promising for the development of discretization procedures for systems in continuous time and space. The aim is to detect stochastic n-point bifurcations (of low order) of the original system via the respective detection in the discretized system. An initial step in this direction such discretizations is done in \cite{Ko13,Ko12}. Rigorous discretization results for stochastic L\'evy (semi-) flows, however, are beyond the scope of the current article and left for future research. However, for the realization of such a project, it is primordial to understand the complexity of such discretized systems, the ground work of which is laid in this article. 

\noindent In this article we answer the following natural questions: 
\begin{enumerate}
 \item[Q.1:] (a) How can a stochastic n-point D-bifurcation be defined rigorously for homogeneous Markovian stochastic (semi-) flows of bijections and mappings? (b) What is the minimal setting in order to define a stochastic n-point D-bifurcation for n-point motions, for instance in situations which do not necessarily come from a stochastic flow? 
 
 \item[Q.2:] (a) What are minimal examples (in the sense of the smallest spaces $\mM$) 
of a stochastic L\'evy semiflow, whose invariant law 
is not determined by its 1- and 2-point motion (opposed to the rigidity of a stochastic Brownian flow)?  
(b) Are there stochastic n-point motions of any level $n$ for any state space of cardinality $m$? 
(c) How can we detect the level of a respective n-point bifurcation 
 given the invariant measure of the flow algorithmically? 
 
\item[Q.3:] (a) How computationally complex can the stochastic n-point D-bifurcations become over a given finite space? That is, how many possible flows of mappings do exist with the same given n-point characteristics? In the context of the rigidity of the laws of stochastic Brownian flows of homeomorphisms: (b) How many linearly independent restrictions are given on every level of n-point motions, and how many of them are necessary to completely determine the law of the respective stochastic flow? 
\end{enumerate}
The article is organized as follows. 
After a review on stochastic flows, Section~\ref{s:notion} provides the setup in terms of a homogeneous Markov n-point system, which generalizes the respective notion of an n-point motion of a homogeneous Markov semiflow to consistent families of transition probabilities in the spirit of \cite{LeJan-Raimond} and \cite{LeJan-Raimond-2}. In Definition \ref{def: n-point bifurcation} of Subsection \ref{ss: def} we introduce the notion of a stochastic n-point D-bifurcation, which completely answers Q.1.

\noindent Section~\ref{s: finite Levy flows} is focussed on the special case of stochastic L\'evy flows over finite sets $\mM$ with $|\mM| = m< \infty$. First we construct two important classes of examples 
of n-point bifurcations in finite spaces with $m$ elements and show stochastic n-point D-bifurcations beyond the well known examples mentioned above (for $n=1,2$) answering Q.2(a) and (b). 
In Subsection \ref{sec: minimal n<m} we present a simple algorithm how to detect the precise level of a stochastic n-point D-bifurcation given the invariant measure of a stochastic L\'evy (semi-) flow. 
We apply this algorithm to two examples of stochastic $3$-point bifurcations, including the minimal one. This provides an answer to Q.2(c). 
The second part of Section \ref{s: finite Levy flows} is devoted to the study of the complexity of stochastic n-point D-bifurcations in Q.3. For stochastic semiflows of mappings we give a complete answer to Q.3(a) in terms of a recursive formula, which is verified for low dimensions by hand and for large values computationally. In case of stochastic flows of bijections, which are the discrete analogue of stochastic Brownian flows of homeomorphisms we conjecture -based on extensive simulations with the data base \cite{OEIS}- a combinatorially interesting, highly nontrivial triangular array of natural numbers $T(m, k)$ described in \cite{Ge90} which quantifies 
to what degree such a stochastic flow is determined by its n-point characteristics ($n\lqq m$). 

Finally, in Subsection \ref{s: continous Levy flows} we embed a stochastic L\'evy flow (and its stochastic $n$-point D-bifurcation) in continuous time and space, in terms of Marcus canonical equation, following the lines of \cite{KPP}.

\bigskip
\begin{footnotesize}
\begin{center}
\begin{tabular}{|l|l|l|l|l|l|}
\hline
\multicolumn{4}{|l|}{Questions: } & \multicolumn{2}{l|}{Answers: } 
\\\hline
\multirow{2}{*}{Q.1} & (a) & \multicolumn{2}{l|}{\multirow{2}{*}{Stochastic n-point D-bifurcation for $n$-point Markov systems}} & 
	\multicolumn{2}{l|}{\multirow{2}{*}{
	% 	Subsection \ref{ss: def}:
	Definition \ref{def: n-point bifurcation}}}
	\\\cline{2-2}
& (b) & \multicolumn{2}{l|}{} & \multicolumn{2}{l|}{}
\\\cline{2-2}
% & (c) & \mult icolumn{2}{c|}{} & \multicolumn{2}{c|}{}
% \\\cline{1-1}
\hline
\multirow{4}{*}{Q.2} & \multirow{2}{*}{(a)} & \multirow{2}{*}{Minimal Examples} & 
for L\'evy (semi-)flows 
& \multicolumn{2}{l|}{%Section \ref{sss: minimal}:
		Example \ref{ex:minimal}}
\\\cline{4-6}
&  &  & 
without a semiflow & \multicolumn{2}{l|}{
% Section \ref{sss: switching}: 
Example \ref{s: nsemiflow}}
\\\cline{2-6}
& (b) & \multicolumn{2}{l|}{Existence of stochastic n-point D-bifurcation} & \multicolumn{2}{l|}{%Section \ref{sss: minimal}:
Lemma \ref{lem-n-point}}
\\\cline{2-6} 
& (c) & \multicolumn{2}{l|}{Level of a stochastic n-point D-bifurcation} & Lemma \ref{cor: projection of invariant measure}  & 
Example \ref{exp: m7}
\\\hline
\multirow{2}{*}{Q.3} & (a) & \multirow{2}{*}{Complexity of the n-pont D-bifurcation} & 
Maps 
		& Recursion \ref{eq: recorrencia} & Pseudocode \ref{algorithm2}
\\ \cline{2-2} \cline{4-6}
& (b) & & 
Bijections & Conjecture \ref{conj:conj} & Pseudocode \ref{algorithm3}
\\\hline
\end{tabular}
\end{center}
\end{footnotesize}

\bigskip
\section{The notion of a stochastic n-point D-bifurcation}\label{s:notion}

We start with some preliminary notation. 
Let $\TT \in \{\NN_0, [0, \infty)\}$ and $\mM$ be a Polish space, that is, a separable, topological space, whose topology is metrizable with respect to a complete metric $\rho$. This includes, $\mM$ being a finite set or the Euclidean space $\RR^d$. Denote by $\cC_\mM = \cC(\mM, \mM)$ the continuous mappings from $\mM\ra \mM$. In case of $\mM$ being discrete and equipped with the discrete topology this coincides with the self-maps $\{\mM \ra \mM\}$. It is equipped with the metric  
 \[
\rho_\infty(f, g) := \sup_{x\in \mM} \rho(f(x), g(x)), \qquad f, g\in \cC_\mM. 
 \]
We denote by $\hH_\mM := \mbox{Homeo}(\mM, \mM)$ the space of homeomorphisms from $\mM\ra \mM$. Note that in case of $\mM$ being discrete with the discrete topology this coincides with the permutations $\mbox{Sym}(\mM) = \{\si: \mM\ra \mM~|~\mbox{ bijective }\}$. It is equipped with the metric 
 \[
 \rho^\infty(f, g) = \rho_\infty(f,g)+ \rho_\infty(f^{-1},g^{-1}), \qquad f, g\in \hH_\mM.  
 \]

 \subsection{The setup and the main notation }
 
We start with the standard setup, see for instance \cite{Ku90, Ku19,LeJan-Raimond}. 
 
\begin{defn} Let $(\Om, \fF, \PP)$ be a complete probability space. 
\begin{enumerate}
 \item A family $(\vp_{s, t})_{\substack{s,t\in \TT\\s \lqq t}}$ of random self-maps on $\mM$ is called a \textbf{Markovian stochastic semiflow of continuous self-maps of $\mM$}, if there is $\nN\in \fF$ with $\PP(\nN) = 0$ such that the following is satisfied:   
 \begin{enumerate}
  \item For all $s, t\in \TT$  with $s\lqq t$ and $\om \in \nN^c$ it follows $\vp_{s, t}(\om) \in \cC_\mM.  $
 \item For all $s, t, u\in \TT$ with $s\lqq u \lqq t$, $\om \in \nN^c$ and $x\in \mM$ it is valid 
 $\vp_{s, t}(\om) = \vp_{u, t}(\om) \circ \vp_{s, u}(\om). $
\item For all $s\in \TT$ and $\om\in \nN^c$ 
 $ \vp_{s, s}(\om) = \id_\mM. $
 \item For all $n\in \NN$, $t_1, \dots, t_n\in \RR$ with $t_1\lqq t_2 \lqq \dots \lqq t_n$ the family of increments $(\vp_{t_1, t_2}, \dots, \vp_{t_{n-1}, t_n}) $ is independent. 
 \item For fixed $s', t'\in \TT$ and $\om\in \nN^c$ the mappings 
 $t \mapsto \vp_{s', t}(\om)$ and $s \mapsto \vp_{s, t'}(\om) $ are c\`adl\`ag (right-continuous with left limits). 
  \end{enumerate}
\item A Markovian stochastic semiflow $(\vp_{s, t})_{\substack{s,t\in \TT\\s \lqq t}}$ of continuous self-maps on $\mM$ is \textbf{homogeneous}~if \\[-5mm]
\begin{enumerate}
 \item[(f)] For all $s,t\in \RR$ with $s\lqq t$ 
 and $h\in \TT$ such that $s+h, t+h\in \TT$ we have 
 $ \vp_{s, t} \stackrel{d}{=} \vp_{s+h, t+h} $.
 \end{enumerate}
A stochastic homogeneous Markovian semiflow $(\vp_{s, t})_{\substack{s,t\in \TT\\s \lqq t}}$ of continuous self-maps of $\mM$ is called a \textbf{stochastic L\'evy semiflow} of continuous self-maps. 
  \item A stochastic L\'evy semiflow $(\vp_{s, t})_{\substack{s,t\in \TT\\s \lqq t}}$ of continuous self-maps of $\mM$ is called a \textbf{stochastic L\'evy flow of homeomorphisms} 
  if the property (a) in item 1. is replaced by 
  \begin{enumerate}
  \item[(a')] For all $s, t\in \TT$  with $s\lqq t$ and $\om \in \nN^c$ it follows 
 $ \vp_{s, t}(\om) \in \hH_\mM.  $
 \end{enumerate}
  \item A stochastic L\'evy semiflow $(\vp_{s, t})_{\substack{s,t\in \TT\\s \lqq t}}$ of continuous self-maps of $\mM$ is called \textbf{stochastic Brownian semiflow of continuous self-maps of $\mM$} 
 or in Kunita's notation \cite{Ku90} a \textbf{Brownian motion with values in $\cC_\mM$} 
 if the property (e) in item 1. is replaced by 
 \begin{enumerate}
  \item[(e')] For fixed $s', t'\in \TT$, $s'\lqq t'$, the mappings 
$ t \mapsto \vp_{s', t}(\om)$ and $s \mapsto \vp_{s, t'}(\om)$ are continuous. 
\end{enumerate} 
 \item A stochastic Brownian semiflow $(\vp_{s, t})_{\substack{s,t\in \TT\\s \lqq t}}$ of continuous self-maps of $\mM$ is called a \textbf{stochastic Brownian flow of homeomorphisms} if the property (a) in item 1. is replaced by (a') of item~3. 
   \item A stochastic Brownian semiflow of continuous self-maps of $\mM$
 is called \textbf{homogeneous} if it satisfies property (f) in item~2. 
 \end{enumerate}
 \end{defn}

\begin{exmp}\label{ex: discrete Levy flow}
{\normalfont 
A stochastic L\'evy semiflow of continuous 
self-maps in discrete time $\TT$ can be written 
as the abstract random walk 
of a random i.i.d. sequence $(\xi_t)_{t\in \TT}$ 
of self-maps in $\mM$. Then, for a positive integer~$t$, the flow $\varphi_t = 
\xi_t(\omega)\circ \ldots \circ \xi_1 (\omega)$. See for instance Arnold \cite{Arnold}, LeJan and Raimond \cite{LeJan-Raimond} 
and references therein.
}
\end{exmp}
\noindent It is one of the main achievements of \cite{Ku90} that stochastic Brownian semiflows of mappings can be characterized as the solution flow of the Fisk-Stratonovich SDE in $\cC_\mM$. 
More general stochastic L\'evy (semi-)flows are found to satisfy the same in case of Marcus canonical equations, see \cite{Ku04, Ku19}. 
\begin{defn}
Given a stochastic homogeneous L\'evy semiflow $\varphi$ of continuous self-maps in $\mM$, $n\in\NN$ and $x = (x_1, \dots, x_n) \in \mM^n$. The respective \textbf{stochastic n-point motion of $\varphi$} is defined by 
\[
\vp_{s,t}(x) := (\vp_{s,t}(x_1), \dots, \vp_{s,t}(x_n)), \qquad t\gqq s\gqq 0.   
\]
\end{defn}
\noindent For a detailed overview we refer the reader to Fujiwara and Kunita \cite{FK85}. Following the lines of the proof of \cite{Ku90}, Theorem 4.2.1, the respective transition probabilities 
\[
P^{(n)}_{s,t}(x, E) := \PP(\vp_{s, t}(x) \in E) 
\]
satisfy the Markov property with respect to the filtration $\fF_{s, t}$ generated by the stochastic L\'evy (semi)-flow~$\vp_{s, t}$.
 \noindent We denote by $\pi^k_{\ell}: \mM^k \ra \mM^{k-1}$, for $1\lqq \ell \lqq k$, 
 the \textbf{projection along the $\ell$-th coordiante} 
% whose image omits the $\ell$-th component 
\[
  \pi^k_ {\ell} (x_1, \dots, x_k) :=  (x_1, \dots,  x_{\ell-1}, x_{\ell+1},  
\dots, x_k) \in \mM^{k-1}.
\]
The following definition turns out to be crucial for the generalization of 
the $n$-point motions of a stochastic flow to situations of merely compatible Markovian families. 
\begin{defn} [Homogeneous n-point Markov System]\label{d:hom Mark Sys} 
Let $\mM$ be a Polish space equipped with its Borel $\sigma$-algebra $\bB(\mM)$ 
and $n\in \NN$ satisfying $n \lqq |\mM|$. 
Consider a family $(P^k)_{1\lqq k\lqq n}$ %indexed by $ k \in \{1,2, \ldots , n\}$ 
of homogeneous transition kernels 
\[
P^k: \TT \times \mM^k \times \bB(\mM^k) \ra [0, 1],  
\]
in the following sense: 
\begin{enumerate}
 \item For any $t\in \TT$, $A\in \bB(\mM^k)$ the map $x\mapsto P^k_t(x, A)$ is measurable. 
 \item For any $t\in \TT$, $x\in \mM^k$ the map $A\mapsto P^k_t(x, A)$ is a probability measure. 
 \item For all $0\lqq s\lqq t$, $x\in \mM^k$ and $A\in \bB(\mM^k)$ the kernel satisfies the Chapman-Kolmogorov equation 
 \[
P^k_t(x, A)= \int_{\mM^k} P^k_{t-s}(z, A)  \ P^k_{s}(x, \mathrm{d} z).   
\]
\end{enumerate}
In addition,  the compatibility of $(P^k)_{1\lqq k\lqq n}$ 
in the sense that all the 
marginals of $P_t^k$ are given by $P_t^{k-1}$, that is, 
\[
P_t^k(x, (\pi^k_ {\ell})^{-1}(A)) = P_t^{k-1}(\pi^k_{\ell} (x), A)\qquad \mbox{ 
for all }t\in \TT\mbox{,  }x\in \mM^k, A\in \bB(\mM^{k-1})\mbox{ and 
}1\lqq \ell \lqq k. 
\]
The pair $(\mM,  (P^k)_ {1\lqq k \lqq n})$ is called a \textbf{homogeneous n-point Markov system}. 
\end{defn}

\begin{exmp}
{\normalfont 
In case of $\mM = \{1, \dots, m\}$, 
a stochastic L\'evy semiflow of mappings $\mM\ra \mM$
and given the respective homogeneous m-point Markov system $(P^k)_{1\lqq k\lqq m}$ 
the law of the flow is uniquely determined by the laws of the respective m-point motions. 
}
\end{exmp}

\begin{exmp}
{\normalfont 
Given a stochastic L\'evy semiflow of mappings over $\mM$, $|\mM| = \{1, \dots, m\}$, $n\lqq m$, the respective family of laws of the $k$-point motions for $1\lqq k\lqq n$ 
forms a homogeneous n-point Markov system.
}
\end{exmp}

\begin{exmp}\label{ex:01noflow}
{\normalfont 
Consider $\mM = \{1, \dots, m\}$ and $n< m$ a homogeneous n-point Markov system $(P^k)_{1\lqq k\lqq n}$ over $\mM$. 
Such laws are not necessarily the distributions of a $n$-point motion 
of a stochastic L\'evy semiflow $\mM \ra \mM$.

\noindent For instance:  A Markov chain in $\mM= \{0,1 \}$ with all transition matrix 
entries equal to $1/2$ and a lifted process in $\mM^2= \{(0,0), (0,1), (1,0), (1,1)\}$ 
with all transition matrix entries equal to $1/4$. 
This system defines a homogeneous 2-point Markov system, 
but a look at the diagonal shows that this dynamics in 
$\mM^2$ is clearly not generated by any semiflow of mappings. 
}
\end{exmp}

The basic object of study in this article is the invariant measure 
of a homogeneous n-point Markov system and its (compatible) projections. 
In context of Definition \ref{d:hom Mark Sys},  for a fixed $n \in \mathbb{N}$, a set $B \in \bB(\mM^n)$ is called $P^n$-{\it invariant}
if $P^n(x,B) = 1$ for all $x \in \mM^n$.  Moreover,
one can define an operator $P^n_\ast\mu$ over the signed measures
on $\bB(\mM^n)$ by
\[
	P^n_\ast \mu (B) = \int_{\mM^n} P^n(y,B)\,d\mu(y).
\]
\begin{defn}
A positive measure $\mu$ on $\bB(\mM^n)$ is called invariant
for a Markov process $X$ with transition probabilities $P^n$ in $\mM^n$ if $P^n_\ast \mu = \mu$.
In addition,  $\mu$ is called ergodic if every invariant set has $\mu$-measure either equal to $0$ or equal to $1$.
\end{defn}

\begin{lem} \label{Prop: projection}
Let $\mu$ be an invariant measure for a Markov process $X$ generated by the 
transition probabilities $P^n$ in $\mM^n$. If 
$\pi^n_k(X)$ is also a Markov process, 
then the induced measure $(\pi^n_k)_*\mu$ is an invariant 
measure for the process $\pi^n_k(X)$ in $\mM^{n-1}$. 
Moreover, if $\mu$ is ergodic then $(\pi^n_k)_*\mu$ is ergodic in $\mM^{k-1}$. 
\end{lem}

\begin{proof} For convenience we drop the superscript $n$ whenever possible. 
It is enough to treat the case $t=1$ which we omit in the sequel. 
Let $P^n(x, A)$ be the family of transition 
probabilities of the process $X$ in $\mM^n$ for $x\in \mM^n$ and  subsets $A\subset 
\mM^n$. The fact that the projection $\pi_k(X)$ generates a Markov process in 
$\mM^{n-1}$ means that the transition probabilities in $\mM^{n-1}$, denoted by  
$P^{n-1}(\pi_k(x), B)$,  is well defined for any $B \subset \mM^{n-1}$ and it is 
given by
\[
P^{n-1}(\pi_k(x), B)= P^n (x, \pi_k^{-1}(B))
\]
for all $x\in \mM^n$. Now, by the induced measure theorem:
\begin{eqnarray*}
 (\pi_k)_* \mu (B)&=& \mu (\pi_k^{-1}(B)) = 
 \int_{\mM^n} P^n(x, \pi_k^{-1}(B))\ d\mu(x)\\ 
&=& \int_{\mM^{n-1}} P^n(\pi_k^{-1}(y), \pi_k^{-1}(B))\ d (p_k)_* \mu(y)= \int_{\mM^{n-1}} P^{n-1}(y, B)\ d (\pi_k)_* \mu(y).
\end{eqnarray*}
The ergodicity follows directly. 
\end{proof}

% \subsubsection{The generalized setup: homogeneous n-point Markov systems}
% \label{ss:setup}

\noindent Note that any homogeneous $n$-point Markov system $(\mM, (P^{k})_{1\lqq k\lqq n})$ can be considered as a compatible extension of $P^1$ from $\mM$ to many possible distributions $P^k$ on $\mM^k$. 
In the sequel we define the stochastic $n$-point motion of a family of a given homogeneous Markov transitions kernel $P$ on $\mM$ (i.e. the characteristics of a $1$-point motion) as any homogeneous $n$-point Markov system $(\mM, (P^{k})_{1\lqq k\lqq n})$ 
with $P^1 = P$ and an additional symmetry condition, which guarantees the indistinguishability of the particles. 

This notion is more general than a stochastic $n$-point motions of a stochastic L\'evy semiflow of mappings, since the latter not necessarily needs to exist as seen in Example~\ref{ex:01noflow}, where the $2$-point motion does not define a semiflow. All such stochastic $n$-point motions define a homogeneous $n$-point Markov system with an additional indistinguishability condition, 
which is also satisfied in case of a stochastic $n$-point motion coming from a homogeneous L\'evy semiflow of mappings.

\begin{defn}\label{def:n-point motion}
% Let $\mM$ be a Polish space, $|\mM|\lqq \infty$, 
% equipped with its Borel sigma algebra $\bB(\mM)$ and denote the time variable by $t\in \TT$, where either $\TT=\mathbb{N}_0$ or $\TT=[0, \infty)$. 
Consider a family $P$ of homogeneous Markov transition kernels 
$P_t (x, A)$, with  $x\in \mM$, $t\in \TT$ and $A \in \bB(\mM)$ %a Borel set in $\mM$ 
which satisfies the classical Chapman-Kolmogorov equation 
\[
P_t(x, A)= \int_{\mM} P_{t-s}(z, A)  \ P_{s}(x, \mathrm{d} z), \qquad 0\lqq s\lqq t.    
\]
For any $n\lqq |\mM|$ a homogeneous $n$-point Markov system $(\mM, (P^k)_{1\lqq k\lqq n})$ which satisfy  
\begin{enumerate}
 \item the \textbf{indistinguishability condition} of the components 
\begin{equation}\label{e:indistignuisheability}
P_t^n ((x_1, \ldots, x_n), B_1\times \dots \times B_n) = P_t^n ((x_{\si(1)}, \ldots, x_{\si(n)}), B_{\si(1)}\times \dots \times B_{\si(n)}), \qquad \si \in S_n, 
\end{equation}
where $(x_1, \ldots, x_n)\in \mM^n$ and $B_1, \dots, B_n$ are Borel sets in $\mM$, and  
 \item the \textbf{extension property} 
\[
P^1 = P, 
\]
\end{enumerate}
is called a \textbf{n-point motion of $P$ (in law)}. 
\end{defn}

\noindent 
\begin{rem}
{\normalfont 
Definition \ref{d:hom Mark Sys} naturally extends to the respective semigroups as follows. 
Given  a homogeneous n-point Markov system 
$(\mM,  (P^k)_ {1\lqq k \lqq n})$ and a continuous bounded function 
$f:\mM^k\rightarrow \mathbb{R}$ the associated Markov semigroup is defined by 
$$
P_t^k f (x)= \int_{\mM^k} f(y) \, P^k_t(x,\mathrm{d} y), \qquad t\in \TT, \ x\in \mM^k.
$$  
The compatibility condition of the family $(P^k)_ {1\lqq k \lqq n} $ in terms of the 
semigroup is equivalent to
\[
P_t^{\ell} f (x_1, \ldots , x_ {\ell})= P_t^k g (y_1, \ldots , y_ k) \qquad 
\mbox{ for all } \ell \lqq k,
\]
whenever $f$ and $g$ are symmetric bounded continuous functions in the sense that 
\[
g (y_1, \ldots , y_ k)= f (y_{i_1}, \ldots , y_{i_ \ell})
\]
for a fixed subset $\{ i_1, \ldots, i_ {\ell}\}\subset \{1, \ldots, k\}$ and 
$(x_1, \ldots , x_ {\ell})= (y_{i_1}, \ldots , y_{i_ {\ell}})$. In the case of 
a compatible family of Feller semigroup, see also e.g. \cite[Def.1.1]{LeJan-Raimond}.
}
\end{rem}

\begin{rem}
{\normalfont
The concept of $n$-point D-bifurcation for a stochastic $n$-point motion of a family of kernels $P$ covers the following three levels of generality: 
\begin{enumerate}
 \item[(1)] The most 
restrictive situation is to assume that the random dynamics comes from a stochastic flow of (measurable) bijections, which are the discrete analogue of stochastic Brownian flows of homeomorphisms. 
 \item[(2)] An intermediate scenario is a stochastic L\'evy semiflow of (measurable) mappings $\varphi$ from $\mM$ to $\mM$, which allows for the  coalescence of particles. As for the bijective case, 
the diagonal is still positively invariant. 
 \item[(3)] However, as shown in Example \ref{ex:01noflow}, in general, the existence of a semiflow cannot be guaranteed. 
 Hence the third perspective consists of stochastic $n$-point motions of a family of homogeneous transition kernel $P$ in the sense of Definition~\ref{def:n-point motion}, independently of any 
flow in $\mM$. As a consequence, sub-diagonals may also no longer be positively invariant, that is, particles can coalesce and ``split''. \\
\end{enumerate}
Our concept of n-point bifurcation applies to all these three cases 
covered by the notion of a n-point motion in the sense of Definition~\ref{def:n-point motion}.  
}
\end{rem}

\bigskip

\subsection{The definition of a stochastic n-point D-bifurcation }
\label{ss: def}
In the theory of deterministic dynamical systems a bifurcation 
occurs when the change of a parameter $\e$ 
of the flow affects the 
support of the invariant 
measure topologically, such as 
for instance splitting into two or more 
disconnected invariant domains. This is well described in classical
dynamical systems, where the precise definition is based on breaking 
local topological equivalences of the flows (see e.g. Katok and Hasselblatt \cite{Katok} 
and the references therein).  For stochastic systems generated by 
It\^o-Stratonovich equations, the bifurcation is mostly considered as the sign change 
of the top Lyapunov exponent, see e.g. L. Arnold \cite{Arnold}, \cite{Arnold-Bifurcation} or recently \cite{La19}.
These two situations have in common the fact that they are observing a breaking 
in the topology of the support of invariant measures, but at different levels: In the 
deterministic  case, the invariant measures are considered in 
$\mM^1$, with trivial extension to $\mM^n$ as the respective product 
measure; in the stochastic case, the sign of the Lyapunov exponents points 
to properties of the invariant measures in $\mM^2$. See the explicit example by Baxendale 
\cite{Baxendale86}, where a bifurcation happens for Brownian motions in the 
torus: i.e. the top Lyapunov exponent change 
the sign but the law of the $1$-point motion is not affected, see also 
\cite{Baxendale91}.

We extend these 1 and 2-point phenomena to n-point motions and 
introduce the following natural generalization 
of a stochastic D-bifurcation for more general 
stochastic flows in Polish spaces. 
Recall that the elements of a family $(A_\e)_{\e>0}$  are topologically equivalent if
for any $\e_1$ and $\e_2$ there is an homeomorphism $h_{\e_1,\e_2}$ such that $A_{\e_1} = h_{\e_1,\e_2}(A_{\e_2})$.
%The notion of topological equivalence among sets comes from the fact that homeomorphic 
%objects can not be distinguished under the topology
%point of view.  This implies there is a natural equivalence relation on these objects which
%is used in the following definition.

\begin{defn}[Stochastic n-point D-bifurcation] \label{def: n-point bifurcation}
Let $\mM$ be a Polish space and $((P^{k, \e})_{1\lqq k\lqq N})_{\e\in I}$
be a family of homogeneous N-point Markov systems in $\mM$ 
indexed by a parameter $\e$ taking values in a real  interval $I$. 
We say, that the family $((P^{k, \e})_{1\lqq k\lqq N})_{\e\in I}$ 
exhibits a \textbf{stochastic n-point D-bifurcation at the bifurcation point $\e_D\in I$} 
for a certain level $n\lqq N$ if it satisfies the following:  
\begin{enumerate}
 \item For any $\bar x\in \mM^k$, $1\lqq k\lqq n$, 
 the mapping $\e \mapsto P^{k, \e}(\bar x, \cdot)$ 
 is continuous with respect to the weak topology on the space of probability measures 
 $\pP(\mM^k)$. 
 \item For any $\e \in I$ 
 there exists an invariant distribution $\mu^{\e}$ 
 with respect to $P^{n, \e}$ on $\mM^{n}$ 
 satisfying the following. 
\begin{enumerate}
 \item For any $\e >\e_D$ the measure $\mu^\e$ is ergodic and 
 all sets of the family $(\supp(\mu^{\e}))_{\e>\e_D}$ 
 are topologically equivalent among each other, but 
\[
\supp (\mu^{\e}) \mbox{ is not topologically equivalent to } \supp(\mu^{\e_D}).
\]
\item For any sequence of projections $\pi^2_{k_2}, \dots, \pi^n_{k_n}$ 
where $k_i \in \{1, \dots, i\}$, $i\in \{2, \dots, n\}$ 
\[
\supp \big((\pi^2_{k_2} \circ \dots \circ \pi^n_{k_n})_* \mu^{\e}\big) \mbox{ is topologically equivalent to }  \supp\big( (\pi^2_{k_2} \circ \dots \circ \pi^n_{k_n})_* \mu^{\e_D}\big).
\]
\end{enumerate}
\end{enumerate}
\end{defn} 
\noindent Examples are given in Subsections \ref{sss: switching}, \ref{sss: minimal} and \ref{sss: arbitrary large order}. 
Since each invariant measure on $\mM^1$ can have many lifts to invariant measures 
in higher levels $\mM^k$, these lifts can exhibit 
more than one stochastic n-point D-bifurcation, 
at different levels $k$. Moreover, the same bifurcation 
on $k$-point motion can have projections into different invariant 
measures on $\mM^1$ (depending on the sequence of projections). 
The precise level n of the n-point D-bifurcation is detected algorithmically in Subsection \ref{sec: minimal n<m}. In Subsection~\ref{ss:degreeoffreedom} the respective numbers of linear restrictions (and its complementary degrees of freedom) are quantified. 

\begin{rem}
{\normalfont 
\begin{enumerate}
 \item As for comparison (following Kunita \cite{Ku90}), we consider 
a homogeneous stochastic Brownian flow $\varphi$ in the group of diffeomorphisms in Euclidean space $\mM = \RR^d$ 
with the infinitesimal mean  
\[
b(x) = \lim_{h\ra 0+} \frac{1}{h} 
\left[\EE\left[\varphi_{h}(x)\right]-x\right], \qquad \forall\; x\in 
\RR^d, 
\]
and the infinitesimal covariance
\[
a(x,y) = \lim_{h\ra 0+} \frac{1}{h} 
\left[\EE\left[\big(\varphi_{h}(x)-x\big) 
\big(\varphi_{h}(y)-y\big)^*\right]\right], \qquad \forall\; x, y\in 
\RR^d. 
\]
Given certain regularity conditions on these parameters (satisfied for instance 
by flows of SDE generated by smooth vector fields with bounded derivatives) the 
law of $\varphi$ in the group of diffeomorphisms is determined by 
$a(x,y)$ and $b(x)$ \cite[Thm. 4.2.5, p. 126]{Ku90}. 
In other words, the law of a (homogeneous) stochastic Brownian semiflow (hence the law of its n-point motion, with $n\gqq 2$) is fully determined by the laws of its 1-point motion and its 2-point motion. 
This result tells us that classical stochastic flows for SDEs driven 
by Brownian motion generically do not furnish the richness of flows 
differing only on higher n-point motion with $n> 2$. 
\item There are several notions of bifurcation in the literature, which are of rather independent nature (for a discussion we refer to \cite{Arnold}, Section 9.1). An alternative definition of a bifurcation would be the sign change of the (leading) Lyapunov exponent. 
In \cite{Baxendale91} the author shows a sign change for the $2$-point motion of the Lyapunov exponent, while leaving the $1$-point characteristics invariant, 
for a stochastic Brownian motion on the torus. While item 1) tells us, that there are at most stochastic $2$-point D-bifurcations for homogeneous stochastic Brownian flows of homeomorphisms, this result indicates that, in fact, there are $2$-point bifurcations for homogeneous stochastic Brownian flows of homeomorphisms. 
\end{enumerate}
}
\end{rem}

\begin{rem}
{\normalfont 
The problem we are addressing here is also related to the recent results \cite{JJost-MKell-CRodrigues} by Jost, Kell and Rodrigues where they study 
conditions under which the transition probabilities (1-point motion) in a 
manifold can be represented by families of random maps. In the same article, 
they consider further conditions for regularity and 
representations by diffeomorphisms. 
For this kind of problem in the context 
of flows which are merely measurable we refer to Kifer 
\cite{Ki86} and Quas \cite{Qu91}.
}
\end{rem}

\begin{rem}
{\normalfont 
The flows we are interested in here are also related to the flow of measurable 
mappings of Le Jan and Raimond \cite{LeJan-Raimond}  (see also \cite{LeJan-Raimond-2})
in the following sense: their flows are constructed from a family of  
Feller compatible semigroups in $C(\mM^n)$, $n\gqq 1$, which preserves the 
diagonal. They are also constructed based on the observation of the 
statistics of the $n$   -point motion, for $n\gqq 1$. 
Problems related to synchronization can also be considered 
in the context for 2- and n-point motions \cite{FGS17, New18}. 
}
\end{rem}

\section{Stochastic n-point D-bifurcations in finite space} \label{s: finite Levy flows}

\subsection{Examples in finite space} 
\label{ss: main example}

In this section we construct different 
examples that exhibit 
a stochastic n-point D-bifurcation of some level $n\gqq 2$. 
The purpose is twofold: We construct classes of examples 
of arbitrarily large cardinality, which are interesting in its own right, 
but also yield the example of a minimal space $\mM$ with $m = |\mM| = 4$ 
announced in the Introduction. The semiflows we are going to address here are time-discrete flow $\varphi_n(x)= \xi^n \circ \ldots \xi^1 (x)$ generated by a composition of a sequence $(\xi^k)_{k\in N}$ of i.i.d. random mappings in $\mM$. Its incremental distribution is precisely the distribution of each random variable $\xi^k, k\in N$. 

 \subsubsection{Minimal example of a stochastic 2-point D-bifurcation (without any semiflow)}\label{sss: switching} 
%   The alternating subgroup $H = A_m$ of $G= S_m$
%  }
Initially we show an example for $m=2$. The novelty here is the notion of a (stochastic) D-bifurcation, even when a semiflow doew not exist. 
% tehre where it does not exist a semiflow. 
% Initially we show an example for $m=2$. The novelty here is, \cblue{ the D-bifurcation}that it does not come from a semiflow. 

\begin{exmp}\label{s: nsemiflow}
{\normalfont As in Example 7, let $(X_n)$ be a Markov chain with state space $\mM = \{0,1\}$ and 1-point transition probability matrix given by
\[
P^1 = \left(\begin{array}{ccc}
            \frac{1}{2} & \frac{1}{2}\\[2pt]
            \frac{1}{2} & \frac{1}{2}\\[2pt]
\end{array}\right),
\]
where $p_{i,j}$ is the transition probability from position $i$ to  position $j$.
Note that $\mu^1 =  \left(\frac{1}{2}, \frac{1}{2}\right)$ is the unique invariant measure of this system. Consider the following family of 2-point motion associated to this system on $\mM^2 = \{(0, 0), (0, 1), (1, 0), (1, 1)\}$ (in lexicographical order) 
which is a consistent  extension from the 1-point motion:
\begin{align*}
P^{2,\e} = \left(\begin{array}{ccccc}
		   \frac{1}{4}+\e & \frac{1}{4}-\e &\frac{1}{4}- \e &\frac{1}{4} + \e\\[2pt]
            \frac{1}{4} & \frac{1}{4} & \frac{1}{4} & \frac{1}{4}\\[2pt]
             \frac{1}{4} & \frac{1}{4} & \frac{1}{4} & \frac{1}{4}\\[2pt]
            \frac{1}{4}+\e &\frac{1}{4}- \e & \frac{1}{4}-\e & \frac{1}{4}+\e\\[2pt]
\end{array}\right)
\end{align*}
 for $\e \in (0,\frac{1}{4})$. The straightforward calculation $\mu^\e P^{2, \e}=\mu^{\e}$ (by Frobenius theorem) shows that the invariant probability measure in $\mM^2$ is given by
\begin{align*}
\mu^\e = \frac{1}{1-2\e}\left(\frac{1}{4}, \,  \frac{1}{4} - \e, \, \frac{1}{4}-\e, \, \frac{1}{4}\right).
\end{align*} 
The system exhibits a stochastic 2-point D-bifurcation in the sense of Definition \ref{def: n-point bifurcation}: if $\e<\e_D = \frac{1}{4}$ then  $\supp(\mu^{\e}) \neq \supp(\mu^{\e_D})$. If a dynamic is generated by a semiflow then it must preserve the diagonal of $\mM^2$. In our context, it means that  $p_{ii,jk}=0$ if $j\neq k$. Note that here $p_{00,01}= p_{00,10}= 1/4 - \varepsilon >0$, for all $\varepsilon \in [0,1/4)$. Hence this example illustrates a D-bifurcation in a dynamics which can not be generated by a semiflow.}

\end{exmp}

\subsubsection{Minimal example of a stochastic 3-point D-bifurcation for a stochastic 
L\'evy flow of bijections}\label{sss: minimal}

Recall that one original motivation of this article is the result by Kunita \cite{Ku90} Theorem~4.2.5 that all stochastic Brownian flows of homeomorphisms in $\RR^d$ under mild conditions on the coefficients are determined by its $1$- and $2$-point characteristics, and hence excludes the existence of $3$-point D-bifurcations. The first result of this subsection shows that the smallest cardinality of the state space $\mM$, which allows for a $3$-point D-bifurcation is $m = |\mM| = 4$.

\begin{lem}\label{lem:minimal}
Over $\mM$ with $m = |\mM| = 4$ 
there exists a stochastic flow of bijections 
which exhibits a $3$-point D-bifurcation, 
while for $m=3$ there are no $3$-point D-bifurcations, 
only $2$ and $1$-point D-bifurcations. 
\end{lem}
\noindent The proof is given by the following example. 
\begin{exmp}[Proof of Lemma~\ref{lem:minimal}]\label{exp: m7}
\label{ex:minimal}
{\normalfont 
Consider $m= 4$ and $G$ the symmetric group $S_4$ over $\mM = \{1,2,3,4\}$. 
Let $H < G$ be the alternating subgroup of $G$ which is given by the even permutations as follows. \\
\begin{center}
\begin{tabular}{|c|c|c|c|}
\hline
No. of transpositions & $H$ & $G\backslash H$ & No. of transpositions \\
\hline
$0$ & $id$ & $(12)$ & $1$ \\
$2$ & $(123) = (12)(23)$   & $(13)$ & $1$ \\
$2$ & $(132) = (13)(32)$   & $(14)$ & $1$ \\
$2$ & $(124) = (12)(24)$   & $(23)$ & $1$ \\
$2$ & $(142) = (14)(42)$   & $(24)$ & $1$ \\
$2$ & $(134) = (13)(34)$   & $(34)$ & $1$ \\
$2$ & $(143) = (14)(43)$   & $(1234) = (12)(23)(34)$ & $3$ \\
$2$ & $(234) = (23)(34)$   & $(1432) = (14)(43)(32)$ & $3$ \\
$2$ & $(243) = (24)(43)$   & $(1324) = (13)(32)(24)$ & $3$ \\
$2$ & $(12)(34)$ &  $(1423) = (14)(42)(23)$ & $3$ \\
$2$ & $(13)(24)$ & $(1243) = (12)(24)(43)$ & $3$ \\
$2$ & $(14)(23)$ & $(1342) = (13)(34)(42)$ & $3$ \\
\hline
\end{tabular}\\ 
\end{center}
We take the uniform distributions $\Delta^H = U(H)$ on $H$ 
and $\Delta^{G\backslash H} = U(G\backslash H)$ and   
consider the discrete flow $\varphi^\e$ associated to the 
increment distribution given by 
\begin{equation} \label{def: perturbed0} 
\Delta_\e := \Delta^H + \e  \left[ \Delta^{G\setminus H} - \Delta^H  \right].
\end{equation}
For $\e>0$ the invariant probability measure is given by $\mu_\e$, which is the uniform distribution 
along the orbits of $G$, while for 
$\e = \e_D:= 0$ we have $\mu_{\e_D}$ is the uniform distribution 
along the orbits of $H$.  
For each pair $(i, j), (k, \ell) \in \mM^2$, $i\neq j$, $k\neq \ell$,  the transition probability $P_{ij,k\ell}^\e$ of the $2$-point motion inherited by the flow which is generated by $\Delta_\e$ 
does not depend on $\e$. In fact, initially note that  there are only two bijections in $G$ 
such that $\mM \setminus \{i, j\} \ra \mM \setminus \{k, \ell\}$, one with an 
even and one with an odd number of transpositions, that is, 
one in $H$ and exactly one in $G\backslash H$. Hence, for $\e\gqq 0$ we have 
\begin{equation}\label{e: pp}
P^\e_{ij, k\ell} = \frac{1-\e}{12} + \frac{\e}{12}= \frac{1}{12}
\end{equation}
where the first summand is the probability of the unique element 
in $H$ occurs and the second summand is the probability that the unique element 
in $G\setminus H$ occurs. Hence the effect of the $\e$-perturbation 
on the transition probabilities in \eqref{def: perturbed0} 
cancels out. 

Therefore, for $\e=0$ and any initial condition outside the subdiagonals 
the trajectories of $\varphi^0$ live almost sure 
in exactly $12 = |H|$ points such that the invariant measure 
satisfies $|\supp(\mu^0)| = |H|=12$. However, for $\e>0$ 
the same argument yields $|\supp(\mu^\e)| = |G|=24$. 
Hence there is an n-point bifurcation of level $n=3$ or $n=4$. 

Note that in general, any initial distribution $\mu$ in $S_n$ and i.i.d. random bijections $\xi^k$ with common distribution $\mu$,  implies that the distribution of the discrete generated flow $\varphi_n = \xi^n \circ \ldots \xi^1$ lies in the subgroup of $S_n$ generated by the support of $ \mu$. Moreover, the  distribution of $\varphi_n$ is the convolution $\mu^{\ast n}= \mu \ast \ldots \mu$, $n$ times. The connection between $\mu $ and the transition probabilities is given by  $\mu\{ \sigma \in S_n : \sigma(i)=k  \mbox{ and } \sigma(j)=l\}= P_{ij,kl}$ for all combinations of $0<i,j,k,l\leq n$.

Hence, in general, to look for  D-bifurcations of level $n=3$ in the  context of (discrete) semiflows of bijections, one has to focus on the support of the action of a proper subgroup of $S_ n$ in $M^2$ and takes an $\epsilon$-perturbation in $\mu$ which enlarge this subgroup, but at the same time, it keeps the same support of invariant measures in $M^2$, or stronger, as above: preserve the transition probabilities in $M^2$, Eq. (\ref{e: pp}). This $\epsilon$-perturbation will enlarge orbits in $M^k$, $k>2$, hence increasing the number of connected components of the invariant measure, i.e. a D-bifurcation happens at level $3<k<n$.

The symmetric group $S_3$ is too small to perform this procedure. It has only 6 elements and 2 proper subsgroups (up to conjugacy) and in neither of them is possible to enlarge the support of $\mu$ without affecting the law, the invariant measure and its support in $M^2$. Hence D-bifurcations with discrete flow of bijections in $\{1,2,3\}$ occurs at level $n=2$, but not higher.
 
In conclusion, $m=4$ represents the minimal number of points 
over which the dynamics of bijections exhibits a stochastic n-point bifurcation 
with $n>2$, which vastly contrasts the rigidity of Brownian flows of diffeomorphisms 
and their invariant measures. 

\noindent This fact is confirmed in Subsection~\ref{sss:dofbijections} quantitatively. 
The upper left corner of Table \ref{tab:dofbijections} reads as follows: 
\begin{center}
\includegraphics[scale=0.6]{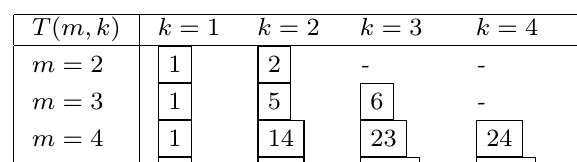}.
\end{center}
Here the number of linear restrictions which imposes the properties of a homogeneous $n$-point 
Markov system is given by the numbers $T(m, k)$, where $n=k$ in this setting. 
In the line for $m=3$, we see that for $n=k = 3>2$, that is, for fixed 2-point characteristics, the $3$ point motion over a space of $3$ elements is (trivially) equal to $6 = m!$. In other words, the matching of the number of variables and linear equations uniquely determines the entire flow.  
However, in the line of $m=4$ we see that the law of the $3$-point motion does not determine the entire flow, since $T(4, 3) = 23 < 24 = 4! = m!$. 
}
\end{exmp}
\noindent The preceding example extends naturally to the following lemma which shows that the bifurcation in higher levels is rather typical.
\begin{lem}\label{lem-n-point}
 Given a finite state space $\mM=\{1,2, \ldots, m\}$, there exist stochastic dynamics with n-point D-bifurcation at level $n=m-1$ or $n=m$ which preserves the transition probabilities of the $(m-2)$-motion.   
\end{lem}
\begin{proof}
We repeat the construction of the preceding example 
with $H= A_m$, being the so-called alternating group, that is 
the subgroup $H< S_m$ with an even number of transpositions. 
There are exactly two elements in $G=S_m$, 
which send any fixed tuple without repetition
$(i_1,\dots, i_{m-2})$ to another fixed tuple without repetition $(j_1, \dots, j_{m-2})$ 
in $\mM^{m-2}$. One of them in $H$ (with an even number of transpositions) 
and one in $G\setminus H$ (with an odd number of transpositions). 
Hence analogously to Example \ref{ex:minimal} the effect of the $\e$-perturbation 
on the $(m-2)$-point transition probabilities cancels out and 
the bifurcation of the invariant measures follows analogously. 
\end{proof}

\subsubsection{A conceptual class of stochastic n-point bifurcations}\label{sss: arbitrary large order}

Following the spirit of the previous subsection we construct 
another family of examples starting with the following basic construction 
extended in the sequel to arbitrary pairs of subgroups $H< G \lneqq S_m$. 

\begin{exmp}\label{exp:m6}
{\normalfont 
 We adopt the notation $f_{i_1i_2i_3i_4i_5i_6}$ for the function 
\[
(1,2,3,4,5,6) \mapsto (i_1,i_2,i_3,i_4,i_5,i_6) 
\]
with $i_1, \ldots, i_m \in \mM$.
The following pairwise notation is convenient for our example:
\[
(1,2,3,4,5,6) = ((1,2),(3,4),(5,6))= (a,b,c) 
\]
and denote by $\bar a, \bar b$ and $\bar c$ the flips of each double entry, 
for example $(\bar a, \bar b, c) = (2,1,4,3,5,6)$. 
Consider the group 
\[
G= \{f_{abc}, f_{\bar abc}, f_{a\bar bc}, f_{ab\bar c}, f_{\bar a\bar bc}, 
f_{\bar 
ab\bar c}, f_{a\bar b\bar c}, f_{\bar a\bar b\bar c}\}
\]
with multiplication given by the composition, and its proper subgroup  
\[
H= \{f_{abc},  f_{\bar a\bar bc}, f_{\bar ab\bar c}, f_{a\bar b\bar c}\}. 
\]
The $6$-point motion of a flow $\vp^0$ generated by composition of  i.i.d. 
bijections with law concentrated on this subgroup, say
\begin{equation}\label{def: unperturbed}
\frac{1}{4} \left[ \delta_{f_{abc}} +  \delta_{f_{\bar a\bar b c}} 
+   \delta_{f_{\bar a b \bar c}} +   \delta_{f_{a\bar b\bar c}} \right],
\end{equation}
has random trajectories with the following property: for each initial condition 
in $\mM^6$ the corresponding random orbit of the process 
is concentrated on at most $4$ points out of $6^6$ possible 
elements of $\mM^6$. For an initial condition which does not belong to any 
subdiagonal (i.e. such that their entries are all different from each other), 
the 
support of the invariant measure is concentrated on exactly $4$ elements. 
On the other hand, the orbits of elements outside any sub-diagonal of the 
$6$-point motion $\vp^\e$ generated by the 
$\e$-perturbation in the law, with $\e >0$,
\begin{equation}\label{eq: perturbed}
\frac{1}{4} 
\left[\delta_{f_{abc}} +   \delta_{f_{\bar a\bar b c}} +   \delta_{f_{\bar a b\bar c}} +   \delta_{f_{a\bar b\bar c}}\right]  
+\frac{\e}{4} \left[ \delta_{f_{\bar abc}} + \delta_{f_{a\bar bc}} + 
\delta_{f_{ab\bar c}} + \delta_{f_{\bar a\bar b\bar c}}
-\delta_{f_{abc}} -   \delta_{f_{\bar a\bar b c}} - \delta_{f_{\bar a b\bar c}} -  \delta_{f_{a\bar b\bar c}} 
\right],
\end{equation}
has invariant measures supported on exactly $8$ elements. Moreover, one easily 
checks by inspection that, due to appropriate cancellations,  the transition 
probability of jumps from a pair of points to any other pair of points does not 
depend on $\e$, i.e the law in $\mM^2$ is constant. The same happens for the law 
in $\mM^1$. 

\noindent The splitting on the number of connected components of the support of the 
invariant measure implies that there must exist an n-point 
bifurcation for $3 \leq n 
\leq 6$. In the next section we shall construct an algorithm to find out  
exactly at 
which level $n$ the bifurcation occurs.
}
\end{exmp}

\paragraph{Extending the construction: } 
For a positive even integer $k$, take the uniform partition of the set 
$\mM=\{1,2, \ldots, m \}$ with $m= k(k+1)$ 
into $(k+1)$ subsets of the form
\begin{align*}
\mM
&=\left\{1,2, \ldots, k \right\} \ \dot\cup\  \left\{(k+1), (k+2), \ldots 2k\right\} \dot\cup  \ \ldots ~ 
\dot\cup \{ k^2+ 1, \ldots, k^2 + k\} = \bigcup_{\ell=1}^{k+1} \{\ell k-k+1, \ldots, \ell k\}.
\end{align*}
For each $\ell \in \{1,2, \ldots, (k+1)\}$, let $ b_\ell$ denote the cyclic 
permutations of the corresponding interval of integers $((\ell k-k+1), \ldots, \ell k)$. 
We consider the Abelian group $G$ of compositions of these cyclic permutations of $\mM$ which 
preserves the subsets of the partition, i.e.:
\[
G= \left\{ b_1^{i_1}\circ b_2^{i_2}\ \circ  \cdots\ \circ b_{k+1}^{i_{k+1}}\ ~|~\ \  
\mbox{ with exponents }\ i_1,i_2,\ldots i_{k+1} \in \{0,1, \ldots k-1\} \    \right\}.
\]
This group has order $|G| = k^{(k+1)}$. Consider the proper subgroup $H<G $ given by
\[
H=\left\{ b_1^{i_1}\ \circ b_2^{i_2}\ \circ \cdots\ \circ b_{k+1}^{i_{k+1}}\ ~|~\ \  
\mbox{ such that  }\ (i_1+i_2 +\ldots ~+i_{k+1}) \mbox{ is even }   \right\},
\]
whose order is $|G|/2$. We take the uniform distribution $\Delta^H$ 
in $H$, i.e. the sum of normalized Dirac measures at each element of $H$. 
Analogously we denote by $\Delta^{G\setminus H}$ the sum of the 
normalized Dirac measures on the $k^{k+1}/2$ elements of its complementary 
set $G\setminus H$. 
With this notation, consider the distribution $\Delta_\e$ in $G$ given by
\begin{equation} \label{def: perturbed} 
\Delta_\e = \Delta^H + \e  \left[ \Delta^{G\setminus H} - \Delta^H  \right] .
\end{equation}
Consider the discrete stochastic flow $\varphi^{\e}$  generated by the composition of i.i.d. random 
elements in $G$ according to the above law. The invariance of the transition probabilities at order less than or equal to $k$ is guaranteed by the following lemma. 
\begin{lem}[There exist stochastic n-point bifurcations of any order]\label{lem: there exist}  
The transition probabilities of the $k$-point motion in $\mM^k$ induced by the 
discrete flow $\varphi^{\e}$ defined above do not depend on $\e>0$. 
Moreover $\varphi^\e$ exhibits a stochastic n-point D-bifurcation for some $n>k$. 
\end{lem}
\begin{proof}
Fix an element $u=(i_1, i_2, \ldots, i_k)$ outside the subdiagonals 
in $\mM^k$. Since the partition of $\mM$ has $(k+1)$ 
subsets and $u$ has $k$ components, 
there exists at least one block $b_\ell$, 
$\ell\in \{1,2, \ldots, (k+1)\}$ whose domain has no intersection with 
 $\{i_1, i_2, \ldots, i_k\}$. 
 Every element 
 $g=b_1^{i_1}\circ b_2^{i_2}\circ  \ldots\circ b_\ell^{i_\ell} \circ \ldots  \circ b_{k+1}^{i_{k+1}}\in G$ acts on $u$ with the following property:  
 \begin{enumerate}
  \item[1)] There are $k/2$ elements 
 $b_1^{i_1} \circ  \ldots \circ b_\ell^{\alpha} \circ \ldots \circ b_{k+1}^{i_{k+1}}\in H$, 
 where $\al \in \{0, \dots, k-1\}$ satisfies the positive parity condition of the exponents, and  
  \item[2)] $k/2$ elements $b_1^{i_1} \circ \ldots \circ b_\ell^{\beta} \circ \ldots \circ b_{k+1}^{i_{k+1}}\in G\setminus H$ for $\beta \in \{0, \dots, k-1\}$, where 
  $\beta$ satisfies the respective negative parity conditions, 
 \end{enumerate}
 such that for any $\alpha$ and $\beta$ given in item 1) and 2) we have 
 \[
 g\cdot u = b_1^{i_1} \circ \ldots\circ b_\ell^{\alpha} \circ \ldots \circ b_{k+1}^{i_{k+1}} \cdot u
 = b_1^{i_1} \circ\ldots \circ b_\ell^{\beta}\circ \ldots \circ b_{k+1}^{i_{k+1}}\cdot u. 
 \]
 Therefore, when one subtracts the probability of the action of 
 elements in $H$ in \eqref{eq: perturbed} the same probability is added to elements in $G\setminus H$ whose 
 action at $u$ is exactly the same. Hence, summing up the independent probabilities 
 that $u$ is sent to any other element in $\mM^k$ does not depend on $\e>0$.
The bifurcation phenomenon of the invariant measures follows analogously as above.  
\end{proof}

\subsection{How to detect the level of a stochastic n-point D-bifurcation}
\label{sec: minimal n<m}

We consider the finite space $\mM = \{1, \dots, m\}$ with $m\gqq 2$. 
The purpose of this section is to answer the following question. 
Given two invariant measures of the $m$-point motion, 
whose projections coincide from a level $k\lqq m$ downwards. 
What is the lowest level $n \in \{k, \dots, m\}$ of 
projections at which they differ? 
The answer of this question yields the level of 
the stochastic n-point D-bifurcation. 
Of course, this can be done along each sequence of projections 
\[
\pi^{k}_{j_k} \circ \dots \circ \pi^{m-1}_{j_{m-1}} \circ \pi^m_{j_m} \mbox{ for }j_\ell \in \{1, \dots, \ell\}, \ell \in \{k, \dots, m\} 
\]
and the detected level will depend on it. 
The minimal level $n> k$ then determines 
the stochastic n-point bifurcation of the entire flow. 

In the sequel we present an algorithm to find 
at which level $n$ the bifurcation happens along a given sequence of projections 
for a given initial condition. 
We apply this algorithm to Example \ref{exp:m6} 
when $k=2$ exhibits a stochastic n-point D-bifurcation 
at levels $n=3$ and $n=5$ for different projections. 
This gives another example of a stochastic 3-point bifurcation 
as motivated in the Introduction, however, 
over a larger than minimal set of points $\mM$ for $m=|\mM| = 6$.

\bigskip
\subsubsection{The projection algorithm}

% For the formulation of the algorithm we need the following some notation. 
For each $1\lqq n\lqq m$ denote by $p^{n}$ be the $m^n \times m^n$ stochastic 
matrix, whose entries are the transition probabilities among the elements of 
$\mM^n$, in the lexicographical order. 
By definition, homogeneous n-point Markov systems 
have transition probabilities which are compatible with projections. 
Hence $p^{n-1}$ can be obtained from the 
projections $\pi^n_r$ for any $r \in \{1, \dots, n\}$ defined 
in Lemma \ref{Prop: projection}. 
More precisely, for  all $1\lqq r \lqq n$ and all 
$(i_1, \dots, i_n), (j_1, \dots, j_n) \in \mM^n\nonumber $ 
we have 
\begin{align}
&p^{n-1}_{\pi^n_r (i_1, \dots  i_n),\pi^n_r (j_1, \dots j_n)} 
= 
\sum_{\ell \in \mM} p^{n}_{(i_1, \dots, i_n),(j_1, \dots,j_{r-1}, \,\ell,\, j_{r+1}, \dots j_n)}.
\label{eq: consistency}
\end{align}
Recall that $\pi^n_r (i_1, \dots,i_r, \dots  i_n) = 
(i_1, \dots, i_{r-1}, i_{r+1}, \dots  i_n) \in \mM^{n-1}$.
% denotes the vector 
% $(i_1, \dots,i_{r-1}, i_{r+1}, \dots  i_n)\in M^{n-1}$. 
This procedure defines a projection $\pi^n_r$ of $P^{n}$ onto $P^{n-1}$ and 
can be expressed algebraically as follows. For each fixed $r$ and $ i_r\in \mM$, there exists a pair of matrices $(R_{n-1}, Q_{n-1})$  
which, according to formula \eqref{eq: consistency}, 
satisfies the equation
\begin{align*}
p^{n-1} = 
R_{n-1} \cdot p^{n}\cdot Q_{n-1},
\end{align*}
where $R_{n-1}$ is a $(m^{n-1} \times m^{n})$-dimensional matrix  with zero  
entries except exactly a unique entry $1$ 
in each row
and $Q_{n-1}$ is an $(m^{n}\times m^{n-1})$-dimensional matrix   
with zero entries except, again,  a unique $1$ in each row. 

The following example illustrates how to use the compatibility, in the sense of Equation (3.6), to reduce transition probabilities in higher order to smaller ones. 
Since there are many choices of projections, there are also many choices of linear combinations of $p_{ij,kl}$. The matrices $R_1$ and $Q1$ below are calculated explicitly using, in each case, the notation presented before for $r$ and $i_r$.
\begin{exmp}
{\normalfont 
For $m=2$, assuming compatibility of the matrix of probability 
transitions $p^{2}$ in $\mM^2$,  we calculate $R_1 \in \RR^{2 \times 4}$ and 
$Q_1 \in \RR^{4 
\times 2}$ for different choices of $(r, i_r)$. 
For $r=1$ and $i_r =1$ we obtain 
\begin{align*}
\left[
\begin{array}{cccc} 
1 & 0 & 0 & 0 \\ 
0 & 1 & 0 & 0 \\  
\end{array}
\right]
\left[
\begin{array}{cccc}
p_{11,11} &  p_{11,12} & p_{11,21} & p_{11,22}\\
p_{12,11} &  p_{12,12} & p_{12,21} & p_{12,22}\\
p_{21,11} &  p_{21,12} & p_{21,21} & p_{21,22}\\
p_{22,11} &  p_{22,12} & p_{22,21} & p_{22,22}\\
\end{array}
\right] 
\left[
\begin{array}{cc} 
1 & 0 \\ 
0 & 1 \\
1 & 0 \\
0 & 1 \\
\end{array}
\right] 
= \left[
\begin{array}{cc}
p_{1,1} & p_{1,2}\\
p_{2,1} & p_{2,2}
 \end{array}
\right].
\end{align*}
Here, $p_{1,1}$ is calculated from point $(1,1)$ as the sum of probability of moving (or not) in the first coordinate: $p_{11,11} +p_{11,21}$; analogously, $p_{1,2}$ is calculated form $(1,1)$ as the sum of probability of moving (or not) in the second coordinate: $p_{11,12} +p_{11,22}$; and so on. For $r=1$ and $i_r =2$: 
\begin{align*}
\left[
\begin{array}{cccc} 
0 & 0 & 1 & 0 \\ 
0 & 0 & 0 & 1 \\  
\end{array}
\right]
\left[
\begin{array}{cccc}
p_{11,11} &  p_{11,12} & p_{11,21} & p_{11,22}\\
p_{12,11} &  p_{12,12} & p_{12,21} & p_{12,22}\\
p_{21,11} &  p_{21,12} & p_{21,21} & p_{21,22}\\
p_{22,11} &  p_{22,12} & p_{22,21} & p_{22,22}\\
\end{array}
\right] 
\left[
\begin{array}{cc} 
1 & 0 \\ 
0 & 1 \\
1 & 0 \\
0 & 1 \\
\end{array}
\right] 
= \left[
\begin{array}{cc}
p_{1,1} & p_{1,2}\\
p_{2,1} & p_{2,2}
 \end{array}
\right].
\end{align*}
Here, $p_{1,1}$ is calculated from point $(2,1)$ as the sum of probability of moving (or not) in the first coordinate: $p_{21,11} +p_{21,21}$; analogously, $p_{1,2}$ is calculated form $(2,1)$ as the sum of probability of moving (or not) in the second coordinate: $p_{21,12} +p_{21,22}$; and so on. Instead, for $r=2$ and $i_r = 1$: 
\begin{align*}
\left[
\begin{array}{cccc} 
1 & 0 & 0 & 0 \\ 
0 & 0 & 1 & 0 \\  
\end{array}
\right]
\left[
\begin{array}{cccc}
p_{11,11} &  p_{11,12} & p_{11,21} & p_{11,22}\\
p_{12,11} &  p_{12,12} & p_{12,21} & p_{12,22}\\
p_{21,11} &  p_{21,12} & p_{21,21} & p_{21,22}\\
p_{22,11} &  p_{22,12} & p_{22,21} & p_{22,22}\\
\end{array}
\right] 
\left[
\begin{array}{cc} 
1 & 0 \\ 
1 & 0 \\
0 & 1 \\
0 & 1 \\
\end{array}
\right] 
= \left[
\begin{array}{cc}
p_{1,1} & p_{1,2}\\
p_{2,1} & p_{2,2}
 \end{array}
\right].
\end{align*}
Note that in all these examples permuting simultaneously 
lines of $R_{1}$ and columns of $Q_1$ leaves the product invariant. 
}
\end{exmp}
\noindent For higher levels, with $m>2$ and $r=i_r=1$ fixed we get the following. 
Thanks to the lexicographical 
order on the entries of the matrices $A^{(n)}$ we have a standard way of 
performing the projections of transition probabilities, using that:
\[
 R_{n-1}=\left[  
     \begin{array}{c|c|c|c}
 
  \iI_{m^{n-1}} &  0 &  \dots  & 0 
  \end{array}
 \right]_ {m^{n-1} \times m^n},
\]
where `$0$' above represents the null $(m^{n-1})$-square matrices, whereas  

\begin{equation}     \label{eq: define Q_{n-1}}
 Q_{n-1}=\left[
     \begin{array}{c|c|c|c}
    \iI_{m^{n-1}} & \iI_{m^{n-1}} & \vdots & \iI_{m^{n-1}} 
  \end{array}
 \right]_{m^{n-1} \times m^{n}}^*.
\end{equation}
Lemma \ref{Prop: projection} implies that given a (left) 
eigenvector $v_n \in \RR^{m^n}$ of $A^{(n)}$, its projection $v_{n-1}$ 
in $\RR^{m^{n-1}}$ is again an eigenvector of 
$A^{(n-1)}$. More precisely we have the following representation. 

\begin{lem} \label{cor: projection of invariant measure}
Given $v_n$ an invariant measure for a compatible Markov chain in the 
product space $\mM^{n}$ represented as a (row)
vector in $\RR^{m^n}$, 
then 
\begin{equation}
v_{n-1} = v_n \ Q_{n-1}  \label{eq: eigenvector projection}
\end{equation}
is an invariant measure in $\mM^{n-1}$ represented as a vector in 
$\RR^{m^{n-1}}$. 
\end{lem}

\begin{proof} Since  formula \eqref{eq: eigenvector 
projection} represents the projection $(\pi^n_1)_*$ in Lemma \ref{Prop: 
projection}. This is straightforwardly, in fact, each column of $Q_{n-1}$ acts on a fixed 
configuration  $(i_2, i_3, \ldots, i_n)\in \mM^{n-1}$, whose 
sum with the first parameter $i_1$ ranging from $1$ to $m$ gives the desired 
    projection.
\end{proof}
\noindent We stress that the matrix $Q_{n-1}$ in formula 
\eqref{eq: eigenvector projection} is not unique and a different choice of $Q_{n-1}$ 
may lead to a different distribution $v_{n-1}$. 
Nevertheless, the choice of $r=1$ and $i_r=1$ 
leads to the simplest version given by \eqref{eq: define Q_{n-1}}.

\bigskip 
\subsubsection{Example \ref{exp:m6} continued: algorithmic detection of a stochastic 3-point bifurcation}\label{ss:detect}

We go back to Example \ref{exp:m6} and apply Lemma \ref{cor: projection of invariant measure} to find at which level $3\lqq n\lqq 6$ the bifurcation occurs given an invariant measure for the
$6$-point motion. 

\begin{lem}
Example \ref{exp:m6} exhibits a stochastic $3$-point D-bifurcations. 
\end{lem}
\noindent The proof is given by the following example. 
\begin{exmp}[Example \ref{exp:m6} continued] 
{\normalfont 
For the sake of notation, we denote by $v^0_\ell$ an invariant 
measure at $\ell$-points for the unperturbed system $\vp^0$ ($\e=0$) 
and by $v^\e_\ell$ an invariant measure of the perturbed system $\vp^\e$, $\e>0$, 
respectively. We represent  these invariant measure as column vector where the points in $\mM^\ell$ are ordered lexicographically  We start with $\ell=6$ and compare $v^0_\ell$ and $v^\e_\ell$, for 
different initial invariant measures. First we consider the 
invariant measures of both systems which contain the 
point $(1,2,3,4,5,6)$. For $\e=0$, the discrete trajectories starting at $(1,2,3,4,5,6)$ pass through points in the set which is exactly the orbit of the action of the subgroup $H$ on this point, hence the trajectories runs recursively (ergodic) in the set (support of invariant measure) $S_H= \{(1,2,3,4,5,6), (1,2,4,3,6,5), (2,1,3,4,6,5), (2,1,4,3,5,6)\}$. And for $\e>0$, the discrete trajectories starting at $(1,2,3,4,5,6)$ pass recursively in the set given by the action of the subgroup $G$ on this point, i.e. in the enlarged  support: $$S_G = S_H \ \bigcup \  \{  (1,2,3,4,6,5), ( 1,2,4,3,6,5),  (2,1,3,4,5,6), (2,1,4,3,6,5)\}.$$
In column representation we obtain 
\[
v^0_6 = \left[
\begin{array}{c} 
 1_{123456} \\ 
 1_{124365} \\
 1_{213465} \\
 1_{214356} \\
\end{array}
\right]_{6^6\times 1} \ \ \ \ \ \ \ \mbox{and} \ \ \ \ \ \ \ \ v^\e_6=\left[
\begin{array}{c} 
 1_{123456} \\ 
 1_{123465} \\
 1_{124356} \\
 1_{124365} \\
 1_{213456} \\
 1_{213465} \\
 1_{214356} \\
 1_{214365}
\end{array}
\right]_{6^6\times 1},
\]
where the symbol $1_{i_1i_2i_3i_4i_5i_6}$ in the column vector means 
that the entry $(i_1,i_2,i_3,i_4,i_5,i_6)$ 
is strictly positive while 
all omitted entries are zero. In addition, we assume as in Example \ref{exp:m6} 
that the distributions $v^0_6$ and $v^\e_6$ are uniform (on their respective support). 

The projections of the invariant measures can easily be performed along the 
first coordinate  $(\pi^\ell_1)_*$, $\ell \in \{2, \dots 6\}$ 
as in Proposition \ref{cor: projection of invariant measure}. 
This means, that according to formula 
($\ref{eq: eigenvector projection}$), one just has to exclude the first 
entry of a nonzero entry $(i_1, \ldots ,i_r)$, $2\lqq r \lqq m $, in 
$v^i_j$, and rearrange, if necessary, in such a way that the order 
in which they appear in the column matrix 
of the reduced level corresponds to the 
lexicographic order again. 
With this method we generate a sequence of vectors, $v_6^0, v^0_5, \dots, v^0_1$, 
which represent each the invariant distributions of a certain level, 
 and $v_6^\e, v^\e_5, \dots , v^\e_1$ 
for the unperturbed and the $\e$-perturbed system accordingly. 
This yields in column vector notation 
\begin{align*}
v_5^0 = Q^T_5 v_6^0 
&= Q^T_{5} \left[
\begin{array}{c} 
 1_{123456} \\ 
 1_{124365} \\
 1_{213465} \\
 1_{214356} \\
\end{array}
\right]_{6^6\times 1} 
= \left[
\begin{array}{c} 
  1_{13465} \\
  1_{14356} \\
  1_{23456} \\ 
  1_{24365} \\
\end{array}
\right]_{6^5\times 1}.
\end{align*}
The complete sequence of projections reads as follows. 

\begin{center}
\begin{tabular}{|l|c|c|c|c|c|c|}
\hline 
%&&&&&&\\
Projections & $v_6^0$ & $v_5^0$ & $v_4^0$ & $v_3^0$ & $v_2^0$ & $v_1^0$ \\
%&&&&&&\\
\hline
&&&&&&\\
Shape of projections&$\left(
\begin{array}{c} 
 1_{123456} \\ 
 1_{124365} \\
 1_{213465} \\
 1_{214356} \\
\end{array}
\right)$
&
$\left(
\begin{array}{c} 
  1_{13465} \\
  1_{14356} \\
  1_{23456} \\ 
  1_{24365} \\
\end{array}
\right)$ 
& 
$\left(
\begin{array}{c} 
 1_{3456} \\
 1_{3465} \\
 1_{4356} \\
 1_{4365} \\
\end{array}
\right)
$ 
&
$
\left(
\begin{array}{c} 
 1_{356} \\
 1_{365} \\
 1_{456} \\
 1_{465}
\end{array}
\right)
$
&
$
\left(
\begin{array}{c} 
 1_{56} \\
 1_{65} 
\end{array}
\right)
$
&
$
\left(
\begin{array}{c} 
% 0 \\ 0 \\ 0  \\0 \\ 
1_{5} \\ 1_{6} \\
\end{array}
\right)
$\\
&&&&&&\\
\hline
Dimension & $6^6\times 1$ & $6^5\times 1$ & $6^4\times 1$ & $6^3\times 1$ & $6^2\times 1$ & $6\times 1$ \\
\hline
\end{tabular}
\end{center}
\noindent For the $\e$-perturbed system we carry out the same algorithm 
and find the invariant measures 
\begin{align*}
v_5^{\e} = Q^T_5 v_6^\e =  
Q^T_{5}
\left[
\begin{array}{c} 
 1_{123456} \\ 
 1_{123465} \\
 1_{124356} \\
 1_{124365} \\
 1_{213456} \\
 1_{213465} \\
 1_{214356} \\
 1_{214365}
\end{array}
\right]_{6^6\times 1}
=
\left[
\begin{array}{c} 
 1_{13456} \\
 1_{13465} \\
 1_{14356} \\
 1_{14365} \\
 1_{23456} \\ 
 1_{23465} \\
 1_{24356} \\
 1_{24365}
\end{array}
\right]_{6^5\times 1} \neq \left[
\begin{array}{c} 
  1_{13465} \\
  1_{14356} \\
  1_{23456} \\ 
  1_{24365} \\
\end{array}
\right]_{6^5\times 1} = v^0_5, 
\end{align*}
\begin{align*}
\mbox{ while } % v_4^\e =  v_4^0 \qquad 
v_j^\e = v_j^0  \ \ \ \ \ \mbox{ for all } \ \ \ \ \ 1\lqq j 
\lqq 
4.
\end{align*}
\noindent This shows that the flow exhibits a stochastic 5-point bifurcation 
for the invariant measure which contains the initial value $(1, 2, 3, 4, 5, 6)$. 

\noindent Taking a different initial invariant measure 
whose support contains the point $(1,2,1,4,1,6)$ 
we have the following projections. 
The sequence of the unperturbed system is given as 
\begin{center}
\begin{tabular}{|l|c|c|c|c|c|c|}
\hline 
Projections & $v_6^0$ & $v_5^0$ & $v_4^0$ & $v_3^0$ & $v_2^0$ & $v_1^0$ \\
\hline
&&&&&&\\
Shape & $\left(
\begin{array}{c} 
 1_{121416} \\ 
 1_{121315} \\
 1_{212425} \\
 1_{212326} 
\end{array}
\right)$
&
$\left(
\begin{array}{c} 
 1_{21416} \\
 1_{21315} \\
 1_{12425} \\ 
 1_{12326}
\end{array}
\right)
$
&
$
\left(
\begin{array}{c} 
 1_{1416} \\
 1_{1315} \\
 1_{2425} \\
 1_{2326} 
\end{array}
\right)
$
&
$\left(
\begin{array}{c} 
 1_{416} \\
 1_{315} \\
 1_{425} \\
 1_{326} 
\end{array}
\right)
$
&
$
\left(
\begin{array}{c} 
 1_{16} \\
 1_{15} \\
 1_{25} \\
 1_{26} 
\end{array}
\right)
$
&
$
\left(
\begin{array}{c} 
% 0 \\ 0 \\ 0  \\0 \\ 
1_{5} \\ 1_{6} \\
\end{array}
\right)
$\\
&&&&&&\\
\hline
\end{tabular}
\end{center}
% \begin{align*}
% v_6^0 = 
% \left[
% \begin{array}{c} 
%  1_{121416} \\ 
%  1_{121315} \\
%  1_{212425} \\
%  1_{212326} 
% \end{array}
% \right]_{6^6\times 1}, 
% v_5^0 
% = \left[
% \begin{array}{c} 
%  1_{21416} \\
%  1_{21315} \\
%  1_{12425} \\ 
%  1_{12326}
% \end{array}
% \right]_{6^5\times 1}, 
% v_4^0 = \left[
% \begin{array}{c} 
%  1_{1416} \\
%  1_{1315} \\
%  1_{2425} \\
%  1_{2326} 
% \end{array}
% \right]_{6^4\times 1},
% \end{align*}
% \begin{align*}
% v_3^0 
% = \left[
% \begin{array}{c} 
%  1_{416} \\
%  1_{315} \\
%  1_{425} \\
%  1_{326} 
% \end{array}
% \right]_{6^3\times 1}, ~
% v_2^0 = \left[
% \begin{array}{c} 
%  1_{16} \\
%  1_{15} \\
%  1_{25} \\
%  1_{26} 
% \end{array}
% \right]_{6^2\times 1}, ~
% v_1^0 = 
% \left[
% \begin{array}{c} 
% % 0 \\ 0 \\ 0  \\0 \\ 
% 1_{5} \\ 1_{6} \\
% \end{array}
% \right]_{6\times 1}.
% \end{align*}
For the $\e$-perturbed system we obtain 
\begin{center}
\begin{tabular}{|l|c|c|c|c|c|c|}
\hline 
Projections & $v_6^\e$ & $v_5^\e$ & $v_4^\e$ & $v_3^\e$ & $v_2^\e$ & $v_1^\e$ \\
\hline
&&&&&&\\
Shape & $\left(
\begin{array}{c} 
 1_{121416} \\
 1_{121415} \\
 1_{121316} \\
 1_{121315} \\
 1_{212426} \\ 
 1_{212425} \\
 1_{212326} \\
 1_{212325} \\
\end{array}
\right)$
&
$
\left(
\begin{array}{c} 
 1_{12426} \\
 1_{12425} \\
 1_{12326} \\
 1_{12325} \\
 1_{21416} \\ 
 1_{21415} \\
 1_{21316} \\
 1_{21315} \\
\end{array}
\right)
$
&
$
\left(
\begin{array}{c} 
 1_{1416} \\
 1_{1415} \\
 1_{1316} \\
 1_{1315} \\
 1_{2426} \\
 1_{2425} \\
 1_{2326} \\
 1_{2325} \\
\end{array}
\right)
$
&
$
\left(
\begin{array}{c} 
 1_{416} \\
 1_{415} \\
 1_{316} \\
 1_{315} \\
 1_{426} \\
 1_{425} \\
 1_{326} \\
 1_{325} \\
\end{array}
\right)
$
&
$\qquad v_2^0\quad$
&
$\qquad v_1^0\quad$\\
\hfill
&&&&&&\\
\hline
Dimension & $6^6\times 1$ & $6^5\times 1$ & $6^4\times 1$ & $6^3\times 1$ & $6^2\times 1$ & $6\times 1$ \\
\hline
\end{tabular} 
\end{center}
% \begin{align*}
% &v_6^\e = \left[
% \begin{array}{c} 
%  1_{121416} \\
%  1_{121415} \\
%  1_{121316} \\
%  1_{121315} \\
%  1_{212426} \\ 
%  1_{212425} \\
%  1_{212326} \\
%  1_{212325} \\
% \end{array}
% \right]_{6^6\times 1}, ~
% v_5^\e =  \left[
% \begin{array}{c} 
%  1_{12426} \\
%  1_{12425} \\
%  1_{12326} \\
%  1_{12325} \\
%  1_{21416} \\ 
%  1_{21415} \\
%  1_{21316} \\
%  1_{21315} \\
% \end{array}
% \right]_{6^5\times 1}, ~
% v_4^\e = 
% \left[
% \begin{array}{c} 
%  1_{1416} \\
%  1_{1415} \\
%  1_{1316} \\
%  1_{1315} \\
%  1_{2426} \\
%  1_{2425} \\
%  1_{2326} \\
%  1_{2325} \\
% \end{array}
% \right]_{6^4\times 1}, ~ 
% v_3^\e = 
% \left[
% \begin{array}{c} 
%  1_{416} \\
%  1_{415} \\
%  1_{316} \\
%  1_{315} \\
%  1_{426} \\
%  1_{425} \\
%  1_{326} \\
%  1_{325} \\
% \end{array}
% \right]_{6^3\times 1}, \\
% &\hfill\\
% &\mbox{ while }v_2^\e = v_2^0, v_1^\e =v_1^0.
%  \end{align*}
\noindent This shows that the flow exhibits in fact 
also a stochastic 3-point D-bifurcation over $m=6$ points,  
in comparison to the minimal Example \ref{ex:minimal} with $m=4$. 
}
\end{exmp}

\bigskip 
\subsection{The complexity of $n$-point D-bifurcations for stochastic L\'evy semiflows}\label{ss:degreeoffreedom} 
%Degrees of freedom to fill up a homogeneous $n$-point Markov systems to a flow}
% Recall that Kunita's theorem states that the two point characteristics of a Brownian flow of homoeomorphisms 
% uniquely determines its law. With this example in mind, it is natural to ask to which degree this behavior fails in case of stochastic L\'evy flows. In contrast to Brownian flows, which cannot be defined in Euclidean space, stochastic L\'evy flows can be defined over a finit space $\mM$. 

In this subsection we study how many linearly independent equations do fixed $k$-point characteristics $P^k$ for all $1\lqq k\lqq n\lqq m$ for some fixed $n$ impose on the laws of the respective flow of self-maps and bijections over $\mM$ (and hence its invariant measure). This combinatorial question is first carried out for the easier case of flows of self-maps $\{\mM\ra \mM\}$ and then in the second case of $\{\mM\ra \mM~|~\mbox{bijective}\}
 = S_m$, $m = |\mM|<\infty$. 
In the first case we prove a recursion formula, which is illustrated numerically. 
In the second case we conjecture - based on numerical experiments - that the respective triangular numbers 
are given by a well-known (complicated) combinatorial quantity introduced in \cite{Ge90}, 
for which to date no closed formula has been found. 
In the appendix we state some explicit formulas for special cases and some asymptotics found there.

\bigskip
\subsubsection{The complexity of $n$-point D-bifurcations for stochastic L\'evy semiflows of random self-maps} \label{sss:dofselfmaps}

% {Degrees of freedom to fill up a homogeneous $n$-point Markov systems to a L\'evy semiflow of random self-maps over a finite space} 

The purpose of this section is to find formulas for the dimensions of the 
vector space of distributions of i.i.d. random self-maps of $\mM= \{1,2, \ldots , m\}$ 
such that the law of the respective flow of random mappings respects the prescribed 
$k$-point characteristics $P^k$ for $0\lqq k\lqq n\lqq m$ for some $n\lqq m$. 

\noindent Recall the notation $f_{i_1 \dots i_m}$ for the function 
$(1, \dots, m) \mapsto 
(i_1, \dots ,i_m)$ for $i_1, \ldots, i_m \in \mM$. The stochastic 
flow of maps $(\vp_{n})_{n\gqq 0}$ in $\mM$ is generated by 
i.i.d. random variables in the space of maps with the following 
discrete probability distribution
\begin{equation}\label{eq: alpha general}
\nu = \sum_{i_1, \dots, i_m =1}^m \alpha_{i_1\dots i_m} \delta_{f_{i_1\dots 
i_m}}, 
\end{equation}
where $\delta_{f_{i_1\dots i_m}}$ is a Dirac measure centered on the mapping 
$f_{i_1\dots i_m}$. The non-negative coefficients $\alpha_{i_1\dots 
i_m} \in \RR$ are ordered lexicographically by the sub-indices. 
We denote 
\[
p_{u_1 \ldots u_k,\ v_1 \ldots v_k} := P_1^k((u_1, \dots, u_k), \{(v_1, \dots, v_k)\}), \qquad (u_1, \dots, u_k),(v_1, \dots, v_k)\in \mM^k. 
\]

\begin{enumerate}
 \item The first linear restriction on the $m^m$ coefficients $(\alpha_{i_1, \dots, i_m})_{(i_1, \dots, i_m)\in \mM^m}$ comes from the fact that they determine 
the distribution of a random variable, hence 
\begin{equation}\label{eq: alpha sum is 1 general}
\sum_{i_1, \dots, i_m=1}^m \alpha_{i_1 \dots i_m} = 1.
\end{equation}
We call this the 0-level restriction for the coefficients. 
\item 
In general, at 
the $k$-level, for a given family of transition probability in $k$-point motion 
$p_{u_1, \ldots, u_k, v_1, \ldots, v_k}$, these characteristics determine linear 
restrictions for the coefficients $\alpha_{i_1, \ldots, i_m}$ given by:
\begin{equation} \label{eq: positions}
 \sum_{(i_1, \dots, i_{m-k})\in \mM^{m-k}} 
\alpha_{(i_1, \dots, i_{m-k}) \lhd \binom{v_1, \dots v_k}{u_1, \dots, u_k}} 
= p_{u_1 
\ldots u_k,\ v_1 \ldots 
v_k},
\end{equation}
where the expression $(i_1, \dots, i_{m-k}) \lhd \binom{v_1, \dots v_k}{u_1, \dots, u_k}$ is the shorthand notation for the following vector 
\[(i_1, \dots, i_{u_1-1}, v_1, i_{u_1+1}, \dots i_{u_2-1}, v_2, i_{u_2+1}, \dots, i_{u_k-1}, v_k, i_{u_k+1}, \dots, i_m).\]
In other words, at position $u_i$ of $(i_1, \dots, i_{m-k})$ the vector $v_i$ is introduced. 
\end{enumerate}
Obviously, the degree of freedom 
(dimension of subspaces which preserve the $k$-point characteristics) is given 
by  $m^m$ minus the number of linearly independent restrictions for the 
coefficients $\alpha_{i_1, \ldots, i_m}$. The following lemma yields a complete recursive description of the number of linearly independent restrictions. 

Given a finite space $\mM = \{1,2, \ldots, m \}$, $1 \lqq k\lqq n\lqq m$ and 
a homogeneous $n$-point Markov system $(P^k)_{1\lqq k \lqq n}$   
we denote by $R^{n, m}_k$ the number of linearly independent restrictions  
imposed simultaneously on the coefficients $(\alpha_{i_1 \dots i_m})$ 
(in the sense of \eqref{eq: positions}) by all $P^\ell$, $\ell\lqq k$, over the alphabet of size $m = |\mM|$. 

\begin{thm}[Recursion formula for the number of restrictions] \label{lem: lin indep restrictions in applications}
In the preceding setting, the triangular numbers 
\[
R^{m, m}_0 \lqq R_1^{m, m}\lqq \dots \lqq R_n^{m, m}
\]
satisfy the following recursion formula:  
For all given $1\lqq k \lqq n\lqq m$ we have 
\begin{align}  \label{eq: recorrencia}
 R^{n, m}_{k} = R^{n, m}_{k-1} + \binom{n}{k}(m^{k} - R^{k, m}_{k-1}), \qquad  R^{n, m}_0 = 1.\\[-5mm]\nonumber
\end{align}
\end{thm}

\begin{rem} 
{\normalfont 
Note that $R^{m, m}_m = m^m$, since the law of the flow of self-maps is uniquely determined, 
that is, the number of variables equals the number of linearly independent equations.\\[-2mm]
}
\end{rem}

\begin{proof}[Proof of Theorem~\ref{lem: lin indep restrictions in applications}:]
% We fix $m\in \NN$ and $1\lqq n\lqq m$.  
We prove by induction over $0\lqq k \lqq n$. 
% that $R^{n, m}_ k$ 
% satisfies the recursive equation \eqref{eq: recorrencia}.
First note that when $k=0$, it means that in formula \eqref{eq: positions} 
the right-hand side is equal to $1$ and on the left-hand side the 
sum is over all $\alpha_{(i_1, \dots, i_m)}$. 
In other words, there are no further restrictions 
other than equation \eqref{eq: alpha sum is 1 general}, that is, $R^{m}_0=1$. 

Assume that the formula holds for $R^{n, m}_{k-1}$, for all $n \in \{k, \ldots ,m\}$. 
The number of restrictions $R^{n, m}_k$ at the $k$-th level depends of the 
characteristics of the level $(k-1)$, i.e. it is a sum of $R^{n, m}_{k-1}$ plus 
some new linearly independent restrictions depending exactly on characteristics at level $k$. 
This justifies the first summand on the right hand side of equation \eqref{eq: recorrencia}. 
It remains to describe these new restrictions depending exactly on characteristics at level $k$. 
By formula \eqref{eq: positions}, considering the projection at level $k$ from 
level $n$ means that there is a subset of positions $\{ \ti u_1, \ldots, \ti u_k\} 
\subseteq \{ u_1, \ldots, u_n\}$ such that for all $v_1, \dots, v_k\in \mM$ we have 
\begin{align} \label{eq: restricao k em n}
& \sum_{
(i_1, \dots, i_{n-k}) \in \mM^{n-k} 
} 
\ti \alpha_{(i_1, \dots, i_{n-k}) \lhd \binom{v_1, \dots, v_k}{\ti u_1, \dots, \ti u_k}} = p_{\ti u_1 
\ldots \ti u_k, \ v_1 \ldots 
v_k}.
\end{align}
The number of possible subsets with $k$ elements $\{\ti u_1, \dots, \ti u_k\}$ out of the ``positions'' $\{ u_1, \ldots, u_n\}$ of the n-point motion yields $\binom{n}{k}$ many ``blocks'' of equations in \eqref{eq: recorrencia} which only vary in $v_1, \dots, v_k\in \mM$. 
However not all of them linearly independent, since $P^k$ inherits dependencies from $P^{k-1}$. 
The number of linearly dependent equations over an alphabet $\mM$ with $|\mM| = m$ at level $k$ inhered form level $k-1$ and below can be expressed in terms of the 
$R^{k, m}_{k-1}$ for which the recursion assumption holds. 
Subtracting $R^{k, m}_{k-1}$ from $m^k$ yields the desired recursion. 
\end{proof}
\bigskip 
\subsubsection{Combinatorial conjecture on the complexity of n-point D-bifurcations for stochastic L\'evy semiflows of random bijections over $m$ points}\label{sss:dofbijections}

In this subsection we conjecture the solution of 
the number of restrictions on the coefficients for flows of 
bijections based on a algorithmic computer experiment and 
\textit{The On-Line Encyclopedia of Integer Sequences} \cite{OEIS}. 
If our conjecture turn out to be true, there is no hope to date to derive 
an analogous recursion formula to Theorem~\ref{lem: lin indep restrictions in applications}. 

\noindent Recall the notation for bijective mappings $f_{i_1 \dots i_m}: \mM 
\ra \mM$. The stochastic L\'evy flow of bijections  $(\vp_{n})_{n\gqq 0}$ 
in $\mM$ is generated by i.i.d. random variables in the space of permutations 
with the following distribution:
\begin{equation}\label{eq: alpha general permutation}
\nu = \sum_{(i_1, \dots, i_m)\in \per(\{1, \dots, m\})} 
\alpha_{i_1\dots i_m} \delta_{f_{i_1\dots i_m}} 
\end{equation}
\begin{enumerate}
\item The first linear restriction on these $m!$ 
coefficients comes from the fact that they determine 
the distribution of a random variable, hence
\begin{equation}\label{eq: alpha sum is 1 general-permutations}
\sum_{(i_1, \dots, i_m) \in \per(\{1, \dots, m\})} \alpha_{i_1 \dots i_m} = 
1.
\end{equation}
We call this the 0-level restriction for the coefficients. 
\item At the level $k$ we have the following. Given a (compatible) family of transition probability in $k$-point motion $p_{u_1 \ldots u_k,\ v_1 \ldots v_k}$, 
we obtain the linear restrictions for the coefficients $\alpha_{i_1, \ldots, i_m}$ given by
\begin{equation} \label{eq: positions permutations}
 \sum_{(i_1, \dots, i_{m-k}) \in \\
\per(\{1, \dots, m \} \setminus \{v_1, \dots, v_k\}) \\ 
} \alpha_{(i_1, \dots, i_{m-k}) \lhd \binom{v_1, \dots, v_k}{u_1, \dots, u_k}} 
= p_{u_1 
\ldots u_k, v_1 \ldots 
v_k},
\end{equation}
where the sum is taken over $(m-k)!$ indices. 
As before, varying the parameters $v_1, \ldots, v_k$ in the expression 
preceding generates a block of $m!/(m-k)!$ equations.
\end{enumerate}
In any level $k$, the diagonal and its complementary are invariant sets for 
the dynamics of random permutations. Moreover, for flows of bijections in 
a finite space, given the sub-maximal $(m-1)$-point transition probabilities, they already determine uniquely the maximal $m$-point transition probabilities, hence they also determine the $m!$ coefficients $\alpha_{i_1 \ldots i_m} $.

Let $u=(u_1, \ldots, u_k)$ and  $v=(v_1, \ldots, v_k)$ be elements in 
$\mM^k$. Since the order of the entries of the elements in $\mM^k$ does not 
matter in a flow, then, if $\sigma$ is a permutation of $k$ 
elements, then by the indistinguishability condition \eqref{e:indistignuisheability}
the transition probabilities satisfy:
\[
p_{u_1 \ldots u_k, v_1 \ldots v_k} = p_{u_{\sigma(1)} \ldots 
u_{\sigma(k)}, v_{\sigma(1)} \ldots v_{\sigma(k)}}.
\]
That is, the entries $(u_1, \dots, u_k)$ can be assumed to be strictly ordered. 

\begin{rem}[A hidden symmetry]
{\normalfont 
Consider now  $u'=(u'_1, \ldots, u'_{(m-k)})$ and $v'=(v'_1, \ldots, 
v'_{(m-k)})$ 
elements in $\mM^{(m-k)}$ such that, as subsets, they  complement $u$ and $v$ 
respectively, 
i.e. $\{u\} \cup \{ u'\} = \{v\} \cup \{ v'\}= \mM$. Then 
\begin{equation} \label{eq: complementaries invariant-permutations}
\sum_{\sigma \in \Delta} p_{u_{1} \ldots 
u_{k}, v_{\sigma(1)} \ldots v_{\sigma(k)}} = \sum_{\xi \in \Delta'} 
p_{u'_{1} \ldots 
u'_{k}, v'_{\xi(1)} \ldots v'_{\xi(k)}}.
\end{equation}
where $\Delta$ are permutations on $k$ elements and $\Delta'$ are permutations 
in $(m-k)$ elements. This is obvious from the observation that in a flow of 
bijections, the whole set $\{ u\}$ is sent to $\{u' \}$ 
(independently of the order), if and only if its complementary $\{v\}$ is sent 
to $\{v'\}$, the complementary of $\{u'\}$. For example:
\[
 p_{1,1} = \sum_{\xi \in  \per(\{2, 3, \ldots, m \})} p_{2  \ldots m\ , \ \xi (2) \xi(3) \ldots \xi(m)}. 
\]
}
\end{rem}

\noindent The arguments in the proofs of Theorem~\ref{lem: lin indep restrictions in 
applications} and of Proposition~\ref{lem: lin indep restrictions permutations} are not easily extensible to higher levels $k$ in the case of bijections. 
This is due to the fact that, in this case, for $k>1$ there are further 
restrictions which  involve crossed equations coming from different blocks of 
equations generated by each fixed $u_1 u_2 \ldots u_k$ in equation 
\eqref{eq: positions permutations}) (in contrast to the previous case of 
arbitrary self-maps). 
Moreover, for level $k\geq m/2$, new restrictions, coming from equation 
\eqref{eq: complementaries invariant-permutations}) which represents 
further dependence on lower levels $(m-k)\leq k$, arises (again different from 
the case of self-maps). Therefore, it looks combinatorially non-trivial 
to control the restrictions coming from different properties with non-empty intersections. 

\noindent However, despite computational issues (see the Appendix) we conjecture with the help of \textit{The On-Line Encyclopedia of Integer Sequences}\cite{OEIS} the numbers $R^{n, m}_k$ for the stochastic L\'evy flows of bijections. 
\begin{conj}\label{conj:conj}
For all $1\lqq k\lqq m$ we have 
\[R^{m, m}_k = m! - T(m, m-k+1),\] 
where the elements of the triangular array $(T(m,k))_{1\leq k\leq m}$ are defined in \cite{OEIS} by 
\begin{center}``Triangle of numbers $T(m,k)$\\ 
$=$ number of permutations of $m$ things with longest increasing subsequence of length $\lqq k$ $(1\lqq k \lqq m)$.''\end{center}
\end{conj}
\noindent The numbers $T(m,k)$ were introduced in \cite{Ge90} and turn out to be highly non-trivial. To date, no recursion or other simple formula is know in the literature for $T(m, k)$. In the Appendix of \cite{Ge90} some asymptotics are derived. On \cite{OEIS1} the triangular array $T(m, k)$ are calculated for $m=1, \dots, 45$. By Remark~\ref{rem:compcomp}, it is plausible, why the algorithmic verification seems to be unfeasible for large values. 

\noindent The following table represents the numbers of $T(m, k)$ as given in \cite{OEIS1}. The $37$ \fbox{boxed} values in the table below have been also calculated by Algorithm~\ref{algorithm3} and coincide with Conjecture~\ref{conj:conj}. \\
\begin{tab}[Algorithmic verification of the triangle numbers $ T(m, k) = m!-R^{m, m}_{m-k+1}$]\label{tab:dofbijections}\hfill
\begin{scriptsize}
\begin{center}
\begin{tabular}{|l|llllllllll|}
\hline
$T(m, k)$&$k=1$&$k=2$& $k=3$ & $k=4$  & $k=5$   & $k=6$   & $k=7$ & $k= 8$ & $k=9$ & $k=10$\\
\hline
$m=2$ & \fbox{$1$}  & \fbox{$2$}   & -     & -      & -       & -       & - & - & -&-\\
$m=3$ &  \fbox{$1$}  & \fbox{$5$}   & \fbox{$6$}   & -      & -       & -       & - & - & -&-\\
$m =4$ &  \fbox{$1$}  & \fbox{$14$}   & \fbox{$23$}  & \fbox{$24$}   & -       & -       & - & - &-&- \\
$m =5$ &  \fbox{$1$}  & \fbox{$42$}  & \fbox{$103$}  & \fbox{$119$}  & \fbox{$120$}   & -       & - & - &-&- \\
$m =6$ &  \fbox{$1$}  & \fbox{$132$}  & \fbox{$513$} & \fbox{$694$}  & \fbox{$719$}  & \fbox{$720$}  & - & - &-&- \\
$m =7$ &  $1$  & \fbox{$429$}  & \fbox{$2761$} & \fbox{$4582$} & \fbox{$5003$} & \fbox{$5039$} & \fbox{$5040$} & - &-&-\\
$m =8$ &  $1$  & $1430$  & $15767$ & $33324$ & \fbox{$39429$} & \fbox{$40270$} & \fbox{$40319$} & \fbox{$40320$}&-&-\\
$m =9$ &  $1$  & $4862$  & $94359$ & $261808$ & $344837$ & \fbox{$361302$} & \fbox{$362815$} & \fbox{$362879$} & \fbox{$362880$}&-\\
$m =10$ &  $1$  & $16796$  & $586590 $ & $2190688$ & $3291590$ & $3587916$ & $3626197$ &\fbox{$3628718$}&\fbox{$3628799$}&\fbox{$3628800$}\\
\hline
\end{tabular}\\
The $37$ \mbox{boxed} values have been verified numerically to coincide with the respective values $T(m, k)$. 
\end{center}
\end{scriptsize}
\end{tab}

\bigskip
\subsection{Embedding of stochastic n-point D-bifurcations on finite space in continuous time and space}\label{s: continous Levy flows}

Stochastic L\'evy semiflows of continuous mappings coming from L\'evy driven SDE are obtained (under rather restrictive conditions, such as finite variation paths) with the help of Marcus canonical equations, the analogue of the Stratonovich differential equations for a class. We refer to Kunita \cite{Ku04}, Section 3, and \cite{Applebaum}, Chapter 6, more details. In discrete time, however, they are given by the random walk representation of Example \ref{ex: discrete Levy flow}. We show in the following example how to embed mutatis mutantis stochastic n-point D-bifurcations in (possibly high-dimensional) Euclidean space. In other words, in the light of Motivation II) of the introduction this is the trivial direction: discrete stochastic n-point bifurcation can be embedded in discrete time and space. The original motivation of approximating stochastic n-point bifurcation of complex systems by its discretizations, however, remains beyond the scope of this article.  

\begin{exmp}
{\normalfont 
In the context of the previous examples with $H < G \leqq S_m$ 
of Example~\ref{lem: there exist} we use the classical representation of the permutations of $S_m$ as elements 
in the orthogonal Lie group $SO(m+1, \RR) \hookrightarrow \mbox{Diff}(\RR^{m+1}, \RR^{m+1})$ 
acting on the first $m$ elements $(e_1, \dots, e_m)$ of the canonical basis of 
the Euclidean space $\RR^{m+1}$. 
To each permutation $f_i\in G$, $i\in \{1, \dots, |G|\}$ we 
associate the unique rotation $U_i$ which sends 
$(e_1, e_ 2, \ldots, e_{m})$ to 
$(e_{f_i(1)}, \ldots , e_{f_i(m)})$. 

By definition, any representation in $SO(m+1, \RR)$ necessarily 
preserves positive orientation. 
However, a generic element $f_i \in G$ may have a negative sign,  
that is, the corresponding $U_i$ has the shape 
\[
\left[ 
\begin{array}{c|c|c|c|c}
&&&& 0\\
e_{f_i(1)} & e_{f_i(2)} & \dots & e_{f_i(m)} & \vdots \\
&&&& 0\\
&&&         & \mbox{sign}(f_i) \\
\end{array}
\right] 
\in \RR^{(m+1) \times (m+1)}.
\]
The last element of the canonical basis  $e_{m +1}$ is sent to $\mbox{sign}(f_i)\cdot e_{m+1}$.
The group $SO(m+1, \mathbb{R})$ is a connected compact Lie group whose Lie algebra 
is the vector space $\mathfrak{so}(m+1)$ of  skew-symmetric matrices. Its compactness guarantees 
that the exponential of elements in the Lie algebra is surjective on $SO(m+1, \mathbb{R})$. 
Hence for each  $U_i \in SO(m+1, \mathbb{R})$, there exists a skew-symmetric matrix $\mathfrak{X}_i$ 
in the Lie algebra of $\mathfrak{so}(m+1, \mathbb{R})$ such that, at time one, the exponential 
satisfies $U_i = \exp {\mathfrak{X}_i}$ for all $i\in \{1, 2, \ldots, |G|\}$. Consider now the following linear Marcus canonical stochastic differential equation \cite{Ma78, Ma81} 
in the generalized Stratonovich sense as in Kurtz, Pardoux, Protter \cite{KPP} 
\begin{eqnarray}\label{e:Mc}
 dx_t & =& \mathfrak{X}_1\, x_t \,  \diamond dZ_t^{1} + \ldots + \mathfrak{X}_{|G|}\,x_t \, \diamond dZ_t^{|G|}
\end{eqnarray}
where the $(Z^{i})_{i=1, \dots, |G|}$ is an i.i.d. family of 
Poisson process (with unitary increment $+1$) of intensity $1 / |G|$ each. A Poisson jump (of height $1$) of $Z^{i}$ in the Marcus equation \eqref{e:Mc} means that the trajectory jumps along the deterministic flow of the corresponding vector field $\mathfrak{X_i}$ for a unitary time, that is, $e^{\mathfrak{X_i}} = U_i$. For further details we refer to \cite{KPP}. As a consequence the process $(x_t)_{t\gqq 0}$ is a compound Poisson process with jump intensity $|G| \cdot 1 /|G| = 1$ 
with a uniform increment distribution on the set $\{U_1, \dots, U_{|G|}\}$. 
Note that Marcus canonical equations are the jump noise equivalent of the Stratonovich 
stochastic integral since it satisfies the Leibniz formula for change of variables 
and hence preserves the flow without additional drift terms, see \cite{KPP}, Proposition 4.2. 

In the sequel we follow the lines of the construction of a stochastic n-point D-bifurcation of Subsection~\ref{sss: arbitrary large order} and keep the respective notation. 
For $H< G$ and $\e =0$ we consider the flow embedding 
in the Marcus sense~\cite{Ma78, Ma81, Applebaum, KPP} 
as above  
\begin{eqnarray}\label{eq:Marcus}
 dx^0_t & =& \mathfrak{X}_1\, x^0_t \,  \diamond dZ_t^{0, 1} + \ldots + \mathfrak{X}_{|H|}\,x^0_t \, \diamond dZ_t^{0, |H|},
\end{eqnarray}
where $(Z^{0, i})_{i\in \{1, \dots , |H|}$ an i.i.d. family of Poisson process 
with intensity $1 / |H|$ each. 
That is, the process $x^0$ is a compound Poisson process with jump intensity $1$ 
and uniform increment distribution on $\{U_1, \dots, U_{|H|}\}$. 
For $\e \in (0, 1]$ we denote by $(Z^{\e, i})_{i\in \{1, \dots, |H|\}}$ an i.i.d. 
family of Poisson processes with intensity $(1-\e) /|H|$ and by 
$(\ti Z^{\e, i})_{i\in \{1, \dots, |G|-|H|\}}$ an i.i.d. 
family of Poisson processes with intensity $\e / (|G|-|H|)$ and consider the following Marcus canonical equation
\begin{align*}
 dx^\e_t 
& = \quad \mathfrak{X}_1 x^\e_t \,  \diamond dZ_t^{\e, 1} + \ldots \quad + \mathfrak{X}_{|H|} x^\e_t \, \diamond dZ_t^{\e, |H|}  + \mathfrak{X}_{|H|+1} x^\e_t \, \diamond d\ti Z_t^{\e, 1} + \ldots  + \mathfrak{X}_{|G|-|H|} x^\e_t \, \diamond d\ti Z_t^{\e, |G|-|H|}. 
\end{align*}
Then $x^\e$ is a compound Poisson process with intensity $|H| \cdot (1-\e)/|H|  + (|G|-|H|) \cdot \e/ (|G|-|H|) =1$ and 
increment distribution 
\begin{equation} \
\Delta_\e := \Delta^H + \e  \left[ \Delta^{G\setminus H} - \Delta^H  \right].
\end{equation}
For $\e=0$ the ergodic invariant measures of the $1$-point motion are given by linear combinations of the Dirac measures centered in the points of the orbit $H.x$, where $x$ is the initial value of the $1$-point motion. Note that for different initial values $x\neq y$ with $H.x \cap H.y = \emptyset$, we have different ergodic measures. 

For $\e>0$ the ergodic invariant measures are obtained analogously with respect to the orbit $G.x$ for each initial value $x$. 
We see from Lemma ~\ref{lem: there exist} 
that each of the ergodic measures of the stochastic flow generated by \eqref{eq:Marcus} with continuous time and space has a stochastic $n$-point D-bifurcation at $\epsilon=0$ for some $n>k$. 
% Hence each ergodic measure associated to an orbit $H.x$ exhibits an $n$-point D-bifurcation. 

In the sequel, we modify the preceding example and obtain a unique invariant measure with smooth density, which exhibits an $n$-point D-bifurcation. 
\medskip
Let $0< r <\frac{\sqrt{2}}{2}$ and take a smooth complete vector field $X$. We assume that the associated deterministic flow is positive invariant on $B(e_i, r)$. 
%  the flow of  $X$ then $\eta^X_t (B(e_i, r_2))\subset B(e_i, r_1)$ for all $t_0< t$ for a certain $0<t_0\in \RR$ and 
  for all $i=1, 2, \ldots, m+1$. Take diffusion coefficients given by a smooth $\sigma: \RR^{m+1} \rightarrow L (\RR^{m+1}; \RR^{m+1})$ such that 
$\sigma (x)$ is nondegenerate (surjective) in the open set $\cup_{i=1}^{m+1} B(e_i, r) $ and has support in the closure $\cup_{i=1}^{m+1} \overline{B(e_i, r)} $.
Consider the system determined by equation
\[
dx_t = X(x_t)\ dt + \sigma(x_t) \circ dW_t, 
\]
where $W = (W_1, \dots, W_{m+1})$ is a $m+1$-dimensional standard Brownian motion. 
% (without jumps). 
By the support theorem (see, Theorem 8.1 in \cite{IW89}), each closed ball $\overline{B(e_i, r)}$ is the support of an ergodic invariant probability measure. Moreover the densities of each of these invariant measures are smooth (inside each ball). 
We keep the previous notation of the vector fields $\mathfrak{X}_i$ and the Poisson processes $(Z_t^{0, 1})_{t\gqq 0}$ and consider for $\e= 0$ the solution of the following equation 
\begin{eqnarray}\label{eq:Marcus2}
 dx^0_t & =& X(x_t^0)\ dt + \sigma(x_t^0) \circ dW_t + \mathfrak{X}_1\, x^0_t \,  \diamond dZ_t^{0, 1} + \ldots + \mathfrak{X}_{|H|}\,x^0_t \, \diamond dZ_t^{0, |H|}
\end{eqnarray}
where $(Z^{0, i})_{i\in \{1, \dots , |H|}$ an i.i.d. family of Poisson process with intensity $1 / |H|$ each. 
For $\e>0$ we add the discontinuous Marcus jump components
\begin{equation}\label{e: M0}
\mathfrak{X}_{|H|+1} x^\e_t \, \diamond d\ti Z_t^{\e, 1} + \ldots  + \mathfrak{X}_{|G|-|H|} x^\e_t \, \diamond d\ti Z_t^{\e, |G|-|H|} 
\end{equation}
to the right-hand side of \eqref{eq:Marcus2} 
as described in the previous paragraph and the superscript $0$ in \eqref{eq:Marcus2} by~$\e$: the L\'evy jumps send the trajectories from one ball $B(e_i, r)$ to another in the right way such that the stochastic D-bifurcation observed in the degenerate case is now reproduced in a setting of non-degenerate smooth densities. 

\bigskip
% That is, the process $x^0$ is a compound Poisson process with jump intensity $1$ 
% and uniform increment distribution on $\{U_1, \dots, U_{|H|}\}$. 
% For $\e \in (0, 1]$ we denote by $(Z^{\e, i})_{i\in \{1, \dots, |H|\}}$ an i.i.d. 
% family of Poisson processes with intensity $(1-\e) /|H|$ and by 
% $(\ti Z^{\e, i})_{i\in \{1, \dots, |G|-|H|\}}$ an i.i.d. 
% family of Poisson processes with intensity $\e / (|G|-|H|)$ and consider 
% \begin{align*}
%  dx^\e_t 
% & = X(x^\e_t)\ dt + \sigma(x^\e_t) \circ dB_t \\
% &\qquad + \mathfrak{X}_1 x^\e_t \,  \diamond dZ_t^{\e, 1} + \ldots \quad + \mathfrak{X}_{|H|} x^\e_t \, \diamond dZ_t^{\e, |H|}  + \mathfrak{X}_{|H|+1} x^\e_t \, \diamond d\ti Z_t^{\e, 1} + \ldots  + \mathfrak{X}_{|G|-|H|} x^\e_t \, \diamond d\ti Z_t^{\e, |G|-|H|}. 
% \end{align*}
% Then $x^\e$ is a compound Poisson process with intensity $|H| \cdot (1-\e)/|H|  + (|G|-|H|) \cdot \e/ (|G|-|H|) =1$. 
}
\end{exmp}

\section*{Acknowledgements: } 
The authors would like to thank Prof. Pedro J. Catuogno for valuable discussions.
The authors acknowledge financial support by DFG within the IRTG 1740: Dynamical Phenomena
in Complex Networks: Fundamentals and Applications. M.A.H. was supported by the FAPA project
``Stochastic dynamics of L\'evy driven systems'' at the School of Sciences of Universidad de los Andes.
M.A.H. also greatfully acknowledges the financial support of Colciencias, the Colombian Administrative
Department of Science, Technology and Innovation for travel support to UNICAMP in July 2019 in the
framework of the Stic Math AmSud2019 project: ``Stochastic analysis of non-Markovian phenomena'',
where this work was completed. M.A.H. also thanks his colleagues Carolina Benedetti and Tristram
Bogart from the Departamento de Matem\'aticas at Universidad de los Andes for their helpful comments
on the project. P.R.R. is partially supported by FAPESP nr. 2015/50122-0, nr. 2020/04426-6 and CNPq
nr. 305212/2019-2. P.H.C. is supported by the project 
Edital DPI/DPG No 03/2020 of Universidade de Brasilia.

\appendix 

\section{Appendix}
\subsection{Pseudocode, examples and calculations of Subsection~\ref{sss:dofselfmaps}}
The following pseudocode algorithm shows how to determine the number of linearly independent restrictions directly. \\
\begin{algorithm}[H]
\caption{Direct algorithmic calculation of number of linearly independent restrictions to complete a homogeneous semiflow of self-maps over $m$-points}\label{algorithm2}
\KwData{homogeneous n-point Markov systems of random self-maps over $\mM=\{1,\ldots,m\}$. }%\ref{eq: positions}}
\KwResult{Number of linearly independent restrictions imposed simultaneously 
on the coefficients $(\alpha_{i_1 \ldots i_m})$ in the sense of~\eqref{eq: positions} imposed by all $P^k$ $k \leqslant n$.}
\quad {\bf For} $k=1$ {\bf to} $n$ {\bf do}\;
\quad\quad \textbf{Set } $A_k =$ \textbf{Matrix}\{coefficients $\alpha_{(i_1, \dots, i_{m-k}) \lhd \binom{v_1, \dots, v_k}{u_1, \dots, u_k}}$ of equations \eqref{eq: positions}\}\\
\quad\quad \textbf{Set }$M_k  =: \left(\begin{array}{c}
            					A_1\\[2pt]
            					A_2\\[2pt]
            					\vdots \\
            					A_k
						\end{array}\right)$ \\
						%(note that each level $k$ we need add new 
						%corresponding equations associated to $P^k$ characteristics)\;
\quad\quad \textbf{Set } $R_k^{n, m} =  \mbox{rank}(M_{k})$\\
\quad {\bf Return } $(R_1^{n, m}, \ldots, R_1^{n, m})$.\\
\end{algorithm}
\bigskip 
\noindent The following table yields the results of the computational illustration of the recursion formula obtained in Theorem~\ref{lem: lin indep restrictions in applications}. 
The unboxed values of the table have been obtained by the recursion formula of Theorem~\ref{lem: lin indep restrictions in applications}. The \fbox{boxed} values have been calculated by Algorithm~\ref{algorithm2} and coincide with the values of the recursion formula of Theorem~\ref{lem: lin indep restrictions in applications}. \\
\begin{tab}[Computational illustration of the recursion formula \eqref{eq: recorrencia}]\label{tab:dofself-maps}\hfill
\begin{scriptsize}
\begin{center}
\begin{tabular}{|l|lllllllllll|}
\hline
$R^{m, m}_n$      & $n=0$ & $n=1$ & $n=2$ & $n=3$  & $n=4$   & $n=5$   & $n=6$ & $n=7$ & $n=8$ & $n=9$ & $n= 10$ \\
\hline
$m =1$ &  \fbox{$1$}  & \fbox{$1^1$}   & -     & -      & -       & -       & - & -& -& -& -\\
$m =2$ &  \fbox{$1$}  & \fbox{$3$}   & \fbox{$2^2$}   & -      & -       & -       & - & -& -& -& -\\
$m =3$ &  \fbox{$1$}  & \fbox{$7$}   & \fbox{$19$}  & \fbox{$3^3$}   & -       & -       & - & -& -& -& -\\
$m =4$ &  \fbox{$1$}  & \fbox{$13$}  & \fbox{$67$}  & \fbox{$175$}  & \fbox{$4^4$}   & -       & -& -& - & -& -\\
$m =5$ &  \fbox{$1$}  & \fbox{$21$}  & \fbox{$181$} & \fbox{$821$}  & \fbox{$2101$}  & \fbox{$5^5$}  & - & -& -& -& -\\%3125
$m =6$ &  \fbox{$1$}  & \fbox{$31$}  & \fbox{$406$} & \fbox{$2906$} & $12281$ & $31031$ & $6^6$ & -& -& -& -\\%46656
$m=7$  &  \fbox{$1$}  & \fbox{$43$}  & \fbox{$799$}  & $8359$ & $53719$ & $217015$ & $543607$ & $7^7$ && -& -\\%823543
$m=8$  &  \fbox{$1$}  & \fbox{$57$}  & $1429$ &$20637$ & $188707$ & $1129899$ & $4424071$ & $11012415$ & $8^8$& -& -\\%16777216
$m=9$ & \fbox{$1$} & $73$ & $2377$ & $45385$ & $561481$ & $4690249$ & $26710345$ & $102207817$ & $253202761$ &$ 9^9$ & -\\%387420489
$m= 10$ & \fbox{$1$} & $91$ & $3736$ & $91216$ & $1469026$ & $16349374$ & $127951984$ & $701908264$ & $2639010709$ 
& $6513215599$ & $10^{10}$\\
\hline
\end{tabular}
\noindent The $31$ values in \fbox{boxes} have been checked numerically. 
% All  values for which the numerical verification was computationally feasible coincide. 
\end{center}
\end{scriptsize}
\end{tab}

\bigskip

 \begin{exmp}[Semiflow of mappings over $m=4$ elements]
 {\normalfont 
 We illustrate the arguments used in the proof of Theorem~\ref{lem: lin indep restrictions in applications}. We consider here flows 
in $\mM=\{1,2,3,4 \}$ with $4^4 = 
256$ coefficients $\alpha_{i_1i_2i_3i_4}$.\\
\textbf{0-point motion:} The number of restrictions is obviously $R^4_0=1$ since it only consists of 
\begin{align*}
\sum_{ijk\ell} \alpha_{ijk\ell} =1. 
\end{align*}
\textbf{1-point motion: }
We have $\binom{4}{1}$ blocks, each block with $m^1 = 4$ new equations:
\begin{align}
\sum_{ijk} \alpha_{uijk} = p_{1,u},\qquad u \in M,\qquad \qquad 
\sum_{ijk} \alpha_{iujk} = p_{2,u},\qquad u \in M,\nonumber\\
\sum_{ijk} \alpha_{ijuk} = p_{3,u},\qquad u \in M,\qquad \qquad 
\sum_{ijk} \alpha_{ijku} = p_{4,u},\qquad u \in M.\nonumber
\end{align}
In each block we have $R^{1, 4}_0=1$ linearly dependent equations which has to be 
subtracted from the total number of equations in the block. Hence, the number of linear  
independent restrictions is given by  
\[
\binom{4}{0}+ \binom{4}{1}[4^1-\binom{1}{0}4^0] = 1 + 4*(4-1) = 13.
\]
\textbf{2-point motion: }
We have $\binom{4}{2} = 6$ blocks, each block with $m^2 = 16$ new equations:
\begin{align*}
\sum_{ij} \alpha_{uvij} = p_{12,uv},\qquad u,v \in M,\qquad \qquad 
\sum_{ij} \alpha_{iuvj} = p_{23,uv},\qquad u,v \in M,\nonumber\\
\sum_{ij} \alpha_{ijuv} = p_{34,uv},\qquad u,v \in M,\qquad \qquad 
\sum_{ij} \alpha_{uivj} = p_{13,uv},\qquad u,v \in M,\nonumber\\
\sum_{ij} \alpha_{iujv} = p_{24,uv},\qquad u,v \in M,\qquad \qquad 
\sum_{ij} \alpha_{uijv} = p_{14,uv},\qquad u,v \in M.\nonumber
\end{align*}
In each block we have $R^{2, 4}_1=7$ linearly dependent equations (obtained 
by putting together the reduction from 2-point motion to 1-point and 0-point motion) which has to be subtracted from the total number of equations of the block. Hence, linearly independent 
restrictions are given by  
\begin{eqnarray*}
\binom{4}{0}[4^0]+ \binom{4}{1}[4^1-\binom{1}{0}4^0] + \binom{4}{2}[4^2 - 
[\binom{2}{0} 4^{0}+ \binom{2}{1} [4^{1} -\binom{1}{0}4^0\Big]\Big] 
% &= &1 + 12 + 6 (16-1 - 2 * 3) \\
% &= & 13+ 6*9 
&= 67.
\end{eqnarray*}
Remaining degrees of freedom $4^4 -67= 256 = 189$.\\
\textbf{3-point motion: }
We have $\binom{4}{3} = 4$ blocks, each block with $m^3 = 64$ new equations:
\begin{align*}
\sum_{i} \alpha_{uvwi} = p_{123,uvw},\qquad u,v, w \in M,\qquad \qquad 
\sum_{i} \alpha_{uviw} = p_{124,uvw},\qquad u,v, w \in M,\nonumber\\
\sum_{j} \alpha_{uivw} = p_{134,uvw},\qquad u,v, w \in M,\qquad \qquad 
\sum_{i} \alpha_{iuvw} = p_{234,uvw},\qquad u,v, w \in M.\nonumber
\end{align*}
In each block we have $R^3_2=37$  linearly dependent equations (obtained 
by putting together the reduction from 3-point 
motion to 2-point, 1-point and 0-point motion) which has to 
be subtracted from the total number of equations of the block. Hence, linearly  
independent restrictions are given by:
\begin{align*}
&\binom{4}{0}[4^0]+\\
&\binom{4}{1}[4^1-\binom{1}{0}[4^{0}]] + \\
&\binom{4}{2}[4^2 -(\binom{2}{2} 4^{0} + \binom{2}{1}[4^{1}-\binom{1}{0}[4^0])] 
+\\
&\binom{4}{1}[4^3- (\binom{3}{3}4^{0} + \binom{3}{2} [4^{1}-\binom{1}{0}[4^0]] 
+ 
\binom{3}{1}[4^{2}-[\binom{2}{0}4^0 +\binom{2}{1}[4^1-\binom{1}{0}[4^0]]]]] 
% = 67 + 4[64-[10+3*9] 
= 175.
\end{align*}
\textbf{4-point motion: }We have a single $\binom{4}{4}=1$ block, with $m^4$ new equations:
\begin{align*}
\alpha_{uvwx} = p_{1234, uvwx}, ,\qquad u,v, w, x \in M,\nonumber\\
\end{align*}
In this single block, we have $R^4_3= 175$  linearly dependent equations 
(obtained by putting together the reduction from 4-point 
motion to 3-point, 2-point, 1-point and 0-point motion) which has to 
be subtracted from the total number of equations of the block. Hence, linearly  
independent restrictions are given by the following equation, which one easily sees that 
has a telescopic cancellation:
\begin{align*}
&{\binom{4}{0}[4^0]}+\\
&{\binom{4}{1}[4^1-\binom{1}{0}[4^{0}]]} + \\
&{\binom{4}{2}[4^2 -[\binom{2}{0} 4^{0} + 
\binom{2}{1}[4^{1}-\binom{1}{0}[4^0]]]} +\\
&{\binom{4}{3}[4^3- [\binom{3}{0}4^{0} + \binom{3}{1} 
[4^{1}-\binom{1}{0}[4^0]] + \binom{3}{2}[4^{2}-[\binom{2}{0}4^0 
+\binom{2}{1}[4^1-\binom{1}{0}[4^0]]]]]}+\\
&\binom{4}{4}[4^4-[{\binom{4}{0}[4^{0}]} + 
{\binom{4}{1}[4^1-\binom{1}{0}[4^0]]} + 
{\binom{4}{2}[4^2- [\binom{2}{0}4^0-\binom{2}{1}[4^1-\binom{1}{0}]}]\\
& +{\binom{4}{3}[4^3- [\binom{3}{0}4^{0} + \binom{3}{1} 
[4^{1}-\binom{1}{0}[4^0]] + \binom{3}{2}[4^{2}-[\binom{2}{0}4^0 
+\binom{2}{1}[4^1-\binom{1}{0}[4^0]]]]]}]] = 4^4.
\end{align*}
The results of this example illustrate the line of $m =4$ in Table~\ref{tab:dofself-maps} for $R^{m, m}_n$, which carries the sequence $(1,13,67,175,256)$. 
}
\end{exmp}
\bigskip 
\subsection{Pseudocode, examples and calculations for Subsection~\ref{sss:dofbijections}}
\begin{rem}[\it The Birkhoff polytopes problem] 
{\normalfont The celebrated 
Birkhoff-von Neumann theorem states that $m\times m$ bi-stochastic matrices lay 
in the 
convex hull generated by  the $m!$ matrices of permutations in $\mM = \{1, 2, 
\ldots,m \}$. For any $m\in \NN$ this convex set is called the Birkhoff polytope $P_m$. There are
 several proofs of this theorem in the literature, for a simple 
and elementary proof we refer to Mirsky \cite{Mirsky}. 
This theory has many interesting application, 
and although already intensely studied, it still offers
some interesting open problems, see Pak \cite{Pak}. 
For instance despite its relevance 
there is no formula for the volume of $P_m$ for higher dimensional 
Birkhoff polytopes $P_m$.  
Only recently an asymptotic formula was obtained by Canfield and McKay 
\cite{Canfield and McKay}. 

In the context of our article 
concerning the random dynamics generated by i.i.d. random mappings  
it means that a stochastic flow in $\mM$ is a flow of permutations 
if and only if the matrices of transition probabilities of 1-point motion 
is not only a stochastic matrix (a matrix whose nonnegative lines entries sum up to $1$), 
but a bi-stochastic matrix (a matrix whose nonnegative lines and column entries sum up to $1$). 
Moreover, in the Birkhoff polytope language, what we are exploring in this 
article is the fact that, in general, except for elements in the wedges of the 
polytope $P_m$, the bi-stochastic matrices has a non-unique representation  as 
a linear combination of the vertices of $P_m$ (in fact, $P_m$ is contained in a 
$(m-1)^2$-dimensional subspace and has $m!$ vertices).
}
\end{rem}
\noindent The flow of bijections induced in the $k$-point level sends each whole fibre 
(component) into a whole fiber. By the Birkhoff-von Neumann theorem  
the matrix of transition probabilities in $k$-point level is 
again bi-stochastic for all $1\leq k \leq m$.
As we have pointed out in the Introduction, in our context, when one deals with 
permutations, it means that one enters in the theory of Birkhoff polytopes, with many open problems. 
For the particular case of $n=2$ we have the following formula. 
\begin{prop}\label{lem: lin indep restrictions permutations}
For a finite space $\mM = \{1,2, \ldots, m \}$,  given 
probability transitions of $1$-point motion, 
the number of linearly independent restrictions 
for the coefficients $(\alpha_{i_1 \dots i_m})$ is given by 
\[
R^{2, m}_1 = (m-1)^2 + 1.
\]
\end{prop}
 \begin{proof}
  The bi-stochastic $m\times m$-matrix of transition probabilities of the 
1-point motion has $2m-1$ redundancies by definition. These redundancies 
correspond to linearly dependent equations of type \eqref{eq: positions 
permutations}) with $k=1$. Hence the restrictions are given by $[m^2 - (2m-1)]$ 
l.i. equations added to the 0-level restriction.
 \end{proof}
As far as our knowledge, in the case of flow of random bijections, the problem of number of 
restrictions on the  coefficients, given the transition probabilities of $k$-point motion, for $k>1$, is still open. 
However the respective numbers of linearly independent restrictions can be calculated numerically. \\

\begin{algorithm}[H]
\caption{Calculation of the linearly independent restrictions of the n-point motion of a L\'evy flow of bijections}\label{algorithm3}
\KwData{stochastic L\'evy flow of bijections over a finite space $\mM=\{1,\ldots,m\}$.  }
\KwResult{Number of linearly independent restrictions imposed simultaneously 
on the coefficients $(\alpha_{i_1 \ldots i_m})$ in the sense of~\eqref{eq: positions permutations} imposed by all $P^k$ $k \leqslant n$.}
\quad {\bf For} $k=1$ {\bf to} $n$ {\bf do}\;
\quad\quad {\bf Set} $A_k =$ \textbf{matrix}\{coefficients $\alpha_{(i_1, \dots, i_{m-k}) \lhd \binom{v_1, \dots, v_k}{u_1, \dots, u_k}}$ of equations \eqref{eq: positions permutations}\}\\
\quad\quad {\bf Set} $M_k  = \left(\begin{array}{c}
            					A_1\\[2pt]
            					A_2\\[2pt]
            					\vdots \\
            					A_k
						\end{array}\right)$ \\
% 						(note that each level $k$ we need add new 
% 						corresponding equations associated to $P^k$ characteristics)\;
\quad\quad {\bf Set} $R^{n, m}_k = \mbox{rank}(M_{k})$\;
\quad {\bf return } $(R^{1, m}_k, \ldots, R^{n, m}_n)$.
\end{algorithm}

\begin{rem}\label{rem:compcomp}
{\normalfont A word about the computational complexity of Algorithms \ref{algorithm2} and \ref{algorithm3}. 
We note that the total amounts of entries of $A^k$ of the 
Algorithm ~\ref{algorithm2} and \ref{algorithm3} grow extremely fast.   
\begin{enumerate} 
\item For the self-maps in Subsection \ref{sss:dofbijections}: 
the total number of variables $(\alpha_{i_1 \ldots i_n})$ (total number of columns) is $n^n$. Moreover for each $P^k$ we have $k!\cdot \binom{n}{k}$  blocks of equations \eqref{eq: positions}
and the amount of (not necessarily linearly independent) 
equation of each block is $n^k$. 
In total the number of equation are $k!\cdot\binom{n}{k}\cdot n^{k}$
consequently the matrix $A^k$ has $k!\cdot \binom{n}{k}\cdot n^{n+k}$ entrances. 
This number shows the difficulty to run computations for large values.

\item For the bijections: the total number of variables $(\alpha_{i_1 \ldots i_n})$ (i.e. columns) is $n!$.  Moreover for each $P^k$
we have $k!\cdot\binom{n}{k}$  blocks of equation and 
the amount of (not necessarily linearly independent) equation of each block is
\[	
n \cdot \ldots \cdot (n-k+1) = k!\cdot \binom{n}{k} = \dfrac{n!}{(n-k)!}.
\]
Then the total number of equation are $\left(\frac{n!}{(n-k)!}\right)^2$
consequently the matrix $A^k$ has $\left(\frac{n!}{(n-k)!}\right)^2 \cdot n!$ 
entrances. In comparison, the numbers in the second case are slightly smaller and hence allow for more values to be verified numerically (31 for self-maps vs. 37 for the bijections). 
\end{enumerate}
}
\end{rem}

\bigskip

\noindent {\bf Acknowledgement:} 
The authors would like to thank the anonymous referee for her / his valuable comments 
which have lead to a significant improvement of the manuscript and the correction of several errors. 
The authors also thank Prof. Pedro J. Catuogno, Prof. C. Benedetti and Prof. T. Bogart for valuable discussions. 

P.H.C. acknowledges support by the project Edital DPI/DPG No 03/2020 of Universidade de Brasilia. 
M.A.H. and P.R.R. acknowledge financial support by DFG within the 
IRTG 1740: {\it Dynamical Phenomena in Complex Networks: Fundamentals and Applications}. 
P.R.R. is partially supported by FAPESP 
nr. 2015/50122-0, nr. 2020/04426-6 and CNPq nr. 305212/2019-2.
M.A.H. was supported by the FAPA project 
``Stochastic dynamics of L\'evy driven systems'' at the School of Sciences of Universidad de los Andes. 
M.A.H. also greatfully acknowledges the financial support of Colciencias, 
the Colombian Administrative Department of Science, Technology and Innovation 
for travel support to UNICAMP in July 2019 
in the framework of the Stic Math AmSud2019 
project: ``Stochastic analysis of non-Markovian phenomena'', 
where this work was completed. 
M.A.H. acknowledges the infrastructure support from the project ``Proyecto de la Convocatoria 2020-2021:
''Stochastic dynamics of systems perturbed with small Markovian noise with applications in
biophysics, climatology and statistics''.

\bibliographystyle{plain}

\end{document}